\documentclass[a4paper, 11pt]{article}
\usepackage{amsmath}
\usepackage{amssymb}
\usepackage{amsthm}
\usepackage{amsfonts}
\usepackage{authblk}
\usepackage{xcolor}
\usepackage[width=460pt, height=650pt]{geometry}
\usepackage{tikz}
\usepackage{soul}
\usepackage[mode=buildnew]{standalone}
\usetikzlibrary{patterns}
\usetikzlibrary{decorations.pathreplacing}
\usetikzlibrary{patterns.meta}
\usetikzlibrary{arrows.meta}
\usetikzlibrary{arrows}
\usepackage[hidelinks]{hyperref}
\usepackage{enumitem}
\usepackage{float}
\usepackage{comment}
\usepackage{microtype}

\providecommand{\keywords}[1]
{
	\small	
	\textbf{\textit{Keywords:}} #1
}

\providecommand{\subjclass}[1]
{
	\small	
	\textbf{\textit{2020 Mathematics Subject Classification:}} #1
}

\newenvironment{Proof}{\noindent \textbf{Proof:}}{\hfill$\square$\\}

\DeclareMathOperator{\N}{\mathbb{N}}
\DeclareMathOperator{\R}{\mathbb{R}}

\DeclareMathOperator{\iii}{\mathfrak{i}}
\DeclareMathOperator{\jjj}{\mathfrak{j}}
\DeclareMathOperator{\kkk}{\mathfrak{k}}
\DeclareMathOperator{\hhh}{\mathfrak{h}}

\DeclareMathOperator{\eps}{\varepsilon}

\DeclareMathOperator{\dH}{\dim_{\mathrm{H}}}
\DeclareMathOperator{\dB}{\dim_{\mathrm{B}}}
\DeclareMathOperator{\dUB}{\overline{\dim}_{\mathrm{B}}}
\DeclareMathOperator{\dLB}{\underline{\dim}_{\mathrm{B}}}

\renewcommand{\phi}{\varphi}

\numberwithin{equation}{section}

\theoremstyle{plain}
\newtheorem{theorem}{Theorem}[section]

\newtheorem{corollary}[theorem]{Corollary}

\newtheorem{lemma}[theorem]{Lemma}

\theoremstyle{remark}
\newtheorem{remark}[theorem]{Remark}

\newtheorem{example}[theorem]{Example}

\theoremstyle{definition}

\tikzdeclarepattern{
	name=mylines,
	parameters={
		\pgfkeysvalueof{/pgf/pattern keys/size},
		\pgfkeysvalueof{/pgf/pattern keys/angle},
		\pgfkeysvalueof{/pgf/pattern keys/line width},
	},
	bounding box={
		(0,-0.5*\pgfkeysvalueof{/pgf/pattern keys/line width}) and
		(\pgfkeysvalueof{/pgf/pattern keys/size},
		0.5*\pgfkeysvalueof{/pgf/pattern keys/line width})},
	tile size={(\pgfkeysvalueof{/pgf/pattern keys/size},
		\pgfkeysvalueof{/pgf/pattern keys/size})},
	tile transformation={rotate=\pgfkeysvalueof{/pgf/pattern keys/angle}},
	defaults={
		size/.initial=5pt,
		angle/.initial=45,
		line width/.initial=.4pt,
	},
	code={
		\draw [line width=\pgfkeysvalueof{/pgf/pattern keys/line width}]
		(0,0) -- (\pgfkeysvalueof{/pgf/pattern keys/size},0);
	},
}

\allowdisplaybreaks

\definecolor{brick}{HTML}{FF0800}

\date{}

\begin{document}
	
	\title{Weakly separated self-affine carpets}
	\author[1]{Bal\'azs B\'ar\'any$^*$}
	\author[1]{Levente David\thanks{BB and LD acknowledges support from the grants NKFI~K142169 and NKFI KKP144059 "Fractal geometry and applications" Research Group.}\thanks{LD received funding from the HUN-REN-BME Stochastics Research Group.}}
	\affil[1]{Department of Stochastics\\ HUN-REN-BME Stochastics Research Group\\ Institute of Mathematics\\  Budapest University of Technology and Economics\\ M\H{u}egyetem rkp. 3., H-
		1111 Budapest, Hungary}
	\affil[ ]{email: barany.balazs@ttk.bme.hu \& ustt1613@gmail.com}
	
	\maketitle

\begin{abstract}
In this paper, we study the Hausdorff and the box-counting dimensions of diagonally aligned self-affine carpets whose projections to the $x$- and $y$-axes satisfy the weak separation condition. In particular, we show that the Hausdorff dimension equals the limit of the Bara\'nski formula, and that the box-counting dimension is the limit of the Feng-Wang formula taken over the $n$-fold compositions of the IFS. We also prove several equivalent formulas for the box-counting dimension, and derive the dimension values for two examples.
\end{abstract}



\section{Introduction}

Let $\mathbb{G}$ be a finite set of contractions over a complete metric space. Such a set is called an iterated function system (IFS). A fundamental result of Hutchinson \cite{Hutchinson1981} asserts that there exists a unique non-empty compact set $\Lambda = \Lambda(\mathbb{G})$ such that $\Lambda=\bigcup_{S\in\mathbb{G}}S(\Lambda)$. We call $\Lambda$ the attractor of $\mathbb{G}$. In $\R^d$, if the maps in $\mathbb{G}$ are affine, we call the IFS $\mathbb{G}$ and its attractor $\Lambda$ self-affine, and in particular, if they are similarities, then we call the IFS and its attractor self-similar.

One of the focal points in the theory of IFSs is the computation of the Hausdorff and the box-counting dimension of the attractor. We now recall the definitions. For $E \subseteq \R^d$, denote the diameter of $E$ (with respect to the usual Euclidean metric) by $|E|$. Let $s \geq 0$, and let us define the $s$-dimensional (outer) Hausdorff measure $\mathcal{H}^s(.)$ by
\begin{align*}
\mathcal{H}^s(E) &:= \lim_{\delta \to 0}\ \inf\bigg\{ \sum_{i \in I}|U_i|^s \ \bigg|\ |U_i| \leq \delta,\ \bigcup_{i\in I}U_i \supseteq E,\ I\text{ is countable}\bigg\}.
\end{align*}
Moreover, we define the Hausdorff dimension of the set $E$ by $\dH(E) := \inf\big\{ s \geq 0 \, \big|\, \mathcal{H}^s(E) = 0 \big\}$.

For a bounded set $E \subset \R^d$, we define the box-counting dimension as
\begin{align*}
\dB(E) :=& \lim_{\delta \to 0^+}\frac{\log N_\delta(E)}{-\log\delta},
\end{align*}
if the limit exists, where $N_\delta(E) := \min\big\{m \in \N\, \big|\, \exists x_1,x_2, \dots ,x_m : E\subseteq \bigcup_{i=1}^m B(x_i, \delta)\big\} $, and $B(x, r)$ denotes the ball of radius $r$ centered at $x$. If the limit does not exist, then the $\liminf$ and the $\limsup$ define the lower- ($\dLB$) and the upper box-counting dimension ($\dUB$). For basic properties of the Hausdorff and box-counting dimension, we refer the reader to Falconer's book \cite{Falconer1990}.

In the case of self-similar IFS on $\R^d$, Hutchinson \cite{Hutchinson1981} showed that $\dUB(\Lambda)\leq\min\{d,s_0\}$, where $s_0$ is called the similarity dimension, and it is the unique solution of the equation $\sum_{S\in\mathbb{G}}r_{S}^{s_0}=1$, where $r_S$ denotes the contraction ratio of the similarity map $S\in\mathbb{G}$. Furthermore, Hutchinson showed that if $\mathbb{G}$ satisfies the \textbf{open set condition} (OSC), that is, there exists a non-empty, open and bounded set $U$ such that $S(U)\subseteq U$ and $S(U)\cap \hat{S}(U)=\emptyset$ for every $S\neq \hat{S}\in\mathbb{G}$, then $\dH(\Lambda)=\dB(\Lambda)=s_0$. The Hausdorff and box-counting dimensions of self-similar sets are equal in general, regardless of their geometric structure, see Falconer \cite{Falconer1997}.

Roughly speaking, the open set condition implies that the overlaps between the images $S(\Lambda)$ (called the cylinder sets) are negligible. The situation becomes more complicated if we allow overlaps. Hochman \cite{Hochman2014, hochman2017selfsimilarsetsoverlapsinverse} showed that if the exponential separation holds, then the Hausdorff dimension equals the minimum of the similarity dimension and $d$. 

Another separation condition is the weak separation condition introduced by Lau and Ngai \cite{LAU199945} and Zerner \cite{MR1343732}. We say that $\mathbb{G}$ satisfies the \textbf{weak separation condition} (WSC) if the identity map is an isolated point of the set
$$
\big\{\hat{S}\circ S^{-1}\, \big|\, \hat{S},S\in\mathbb{G}^*\big\},
$$
where $\mathbb{G}^*=\{S_1\circ\cdots\circ S_n\, |\,S_1,\ldots,S_n\in\mathbb{G},\, n\in\N\}$ is the semigroup induced by the maps in $\mathbb{G}$. The WSC allows exact overlaps, while the non-identical cylinders remain relatively well separated. This property enables the computation of the Hausdorff dimension of the attractor by taking the limit of the similarity dimensions associated with higher iterates of the IFS, once exact overlaps have been removed.

An equivalent characterisation of the WSC is the existence of a constant $C$ such that every ball of a given radius intersects at most $C$ cylinders of comparable size. This property ensures that the natural cover used to compute the Hausdorff dimension, namely, the cover consisting of images of balls under higher iterates of the defining maps, remains asymptotically optimal. The same conclusion continues to hold if a subpolynomial bound replaces the constant $C$ in terms of radii. This weaker requirement is called the asymptotic weak separation condition and was introduced in Feng \cite{MR2322179}.

Diagonally aligned carpets are among the simplest examples of self-affine sets that are non-self-similar. One of the first studies of such sets is due to Bedford \cite{bedford1984crinkly} and McMullen \cite{McMullen}, who independently calculated the Hausdorff and box-counting dimensions of self-affine sets with contractions mapping the unit square to rectangles in a homogeneous rectangular grid. In contrast to the result of Falconer \cite{Falconer1997}, we observe that in this setting the equality of the Hausdorff and box-counting dimensions is not typical, and can be characterised with a simple geometric condition.

Their construction has many constraints, which led to generalisations along various avenues. Lalley and Gatzouras \cite{MR1183358} considered planar diagonally aligned self-affine carpets with a specific column structure and a strict ordering of the contraction rates along the coordinate axes. Later, Bara\'nski \cite{MR2298824} generalised the notion to allow cases in which some functions contract more along the first coordinate axis, and others contract more along the second. In both cases, the Hausdorff and box-counting dimensions were determined. Feng and Hu \cite{MR2560042} proved a formula for the Hausdorff dimension of self-affine ergodic measures supported on diagonally aligned self-affine sets, which establishes relations among projected entropies, Lyapunov exponents, and the dimensions of projections onto the coordinate axes. This result is the base expression for the various formulas for the Hausdorff dimension of the attractor. In particular, Lalley and Gatzouras \cite{MR1183358} and Bara\'nski \cite{MR2298824} proved that the Hausdorff dimension of the attractor is the maximum of the Feng-Hu formula taken over the probability distributions on $\mathbb{G}$ with Lyapunov exponents corresponding to the contraction ratios along the $x$- and $y$-axes.

In a more general setup, namely, without any grid-like structure or order between the contraction ratios, Feng and Wang \cite{MR2128947} determined the box-counting dimension under the rectangular open set condition for self-affine carpets.

Several articles have recently been published on the dimension theory of overlapping self-affine carpets. Fraser and Shmerkin \cite{MR3570020} modified the construction of Bedford-McMullen by considering typical translations of the columns formed by the maps, and Pardo-Sim\'on \cite{simon2017dimensionsoverlappinggeneralizationbaranski} considered overlapping Bara\'nski carpets by taking random translations of the rows and columns formed by the maps. In particular, they assumed that the projections of the columns (and rows) satisfy the exponential separation condition for self-similar IFSs. Recently, Rapaport \cite{rapaport2023dimensiondiagonalselfaffinesets} and Feng \cite{feng2025dimensiondiagonalselfaffinemeasures} considered general overlapping diagonal self-affine sets under the exponential separation condition of the iterated function systems induced by the coordinate projections.

We wish to continue the study of overlapping self-affine carpets on the plane by introducing a variation of the asymptotic weak separation condition for our construction. Our standing assumption is that the coordinate projections (as self-similar systems) satisfy the asymptotic weak separation condition. Under this assumption, we provide a formula for the Hausdorff and the box-counting dimension of the attractor.

Our results can also be considered a generalisation of the planar case of the results of He, Lau and Rao \cite{He2003}, who considered systems with homogeneous linear parts that are inverses of expanding integer-coefficient matrices with translation vectors of integer coordinates. Such systems are strongly related to sofic self-affine fractals, see for example Kenyon and Peres \cite{Kenyon1996, MR1389626}, and Alibabaei \cite{alibabaei2024exacthausdorffdimensionsofic}.

\subsection{Setup}
We now present our main assumptions and findings. First, we introduce the assumptions on the systems under consideration.

It is standard that any affine map $S:\R^2 \to \R^2$ is of the form $S(x) := A_S\,x+t_S$, where $A_S$ is a $2\times2$ real matrix and $t_S\in \R^2$.

\begin{enumerate}[label = $\circ$, ref= self-affine carpet]
	\item \label{A1} Let $\mathbb{G}$ be a finite collection of maps of the form
\begin{equation}\label{eq:mapform}
S(x,y)=\left(r_{S,1}x+t_{S,1},r_{S,2}y+t_{S,2}\right),
\end{equation}
such that $|r_{S,\ell}|\in(0,1)$ for every $S\in\mathbb{G}$ and $\ell\in\{1,2\}$. That is, $A_S$ is diagonal for any $S \in \mathbb{G}$. We call such an IFS and its attractor a diagonal or diagonally aligned self-affine set, or simply \textbf{self-affine carpet}.
\end{enumerate}

Denote the orthogonal projections to the main coordinate axes by $\mathrm{p}_1(x,y)=x$ and $\mathrm{p}_2(x,y)=y$, respectively. These projections are called the principal \textbf{projections}. It is easy to see that the orthogonal projections of the maps in $\mathbb{G}$ form self-similar IFSs on $\R$ with attractor $\mathrm{p}_\ell\Lambda$, $\ell \in \{1,2\}$. Indeed, for any map $S$ of the form in \eqref{eq:mapform}, $\mathrm{p}_1\circ S(x,y)$ is independent of $y$ and so we can define $\mathrm{p}_1S(x):=\mathrm{p}_1\circ S(x,y)=r_{S,1}x+t_{S,1}$ (respectively for $\mathrm{p}_2S$). Let us denote this IFS by $\mathrm{p}_\ell\mathbb{G}:=\{\mathrm{p}_{\ell}S: S\in\mathbb{G}\}$. 

Let us note that it may happen that $\mathrm{p}_\ell S = \mathrm{p}_\ell\hat{S}$, even though $S\neq\hat{S}\in\mathbb{G}$. To avoid redundancies in projections and higher iterates, we always consider $\mathrm{p}_\ell\mathbb{G}$ as a set.

\begin{enumerate}[label = $\circ$, ref= CC]
\item\label{A2} We assume that $\mathrm{conv}(\mathrm{p}_1\Lambda) \times \mathrm{conv}(\mathrm{p}_2\Lambda) = [0,1]^2$, where $\mathrm{conv}(\cdot)$ refers to the convex hull of a set.
\end{enumerate}

The above assumption is purely technical; it excludes cases in which the attractor lies on a horizontal or vertical line. In those cases, however, the IFS is self-similar and has a well-developed but distinct theory.

Concerning self-affine carpets, we may consider the following additional structural conditions:

\begin{enumerate}[label = $\circ$, ref = homogeneous]
\item \label{B1} We say that an IFS has \textbf{homogeneous contractions} if there exist positive reals $r_1, r_2$ such that $r_{S,1}=r_1$ and $r_{S,2}=r_2$ for every $S\in\mathbb{G}$. Without loss of generality, we always assume $r_1\leq r_2$ in this case.
\end{enumerate}

\begin{enumerate}[label = $\circ$, ref = orientation preserving]
\item\label{B2} We say that the IFS are \textbf{orientation preserving} if both $r_{S,1}$ and $r_{S,2}$ are positive for every $S\in\mathbb{G}$.
\end{enumerate}

\begin{enumerate}[label = $\circ$, ref = ROSC]
\item\label{C1} We say that $\mathbb{G}$ satisfies the \textbf{rectangular open set condition} (ROSC) if $S((0,1)^2) \cap \hat{S}((0,1)^2) = \emptyset$ for any $S\neq \hat{S} \in \mathbb{G}$.
\end{enumerate}

\begin{enumerate}[label = $\circ$, ref = coordinate ordering]
\item\label{G1} We say that $\mathbb{G}$ satisfies \textbf{coordinate ordering} if $|r_{S,1}| \leq |r_{S,2}|$ for any $S\in\mathbb{G}$.
\end{enumerate}

Note that the previously studied cases of Bedford \cite{bedford1984crinkly} and McMullen \cite{McMullen}, Lalley and Gatzouras \cite{MR1183358}, and Bara\'nski \cite{MR2298824} can be described by these structural assumptions. For example, a \ref{B1}, \ref{B2}, \ref{A1} with $r_i^{-1}$ a positive integer, satisfying \ref{C1} and the OSC for both projections, essentially describes Bedford-McMullen carpets.

An \ref{B2} \ref{A1} satisfying the \ref{C1}, \ref{G1} and the OSC for the $y$-projection is a Lalley and Gatzouras~\cite{MR1183358} carpet; and a \ref{B2} \ref{A1} satisfying the \ref{C1} and the OSC for both projections corresponds to Bara\'nski \cite{MR2298824} carpets. We will call these setups the well-separated cases.

The \ref{B2} assumption is inessential. In the case of Lalley and Gatzouras \cite{MR1183358}, their arguments would work if rows of functions are allowed to switch orientation simultaneously along their $y$-axes, while switching along the $x$-axes appears insignificant. In the case of Bara\'nski carpets \cite{MR2298824}, a similar extension would need more attention.

\subsection{Main results}

Here we state a streamlined, concise form of our theorems, using the WSC and \ref{C1}. We note that in their respective sections, these theorems take a much more overgrown form, using separation conditions that rely on the asymptotic weak separation condition (\ref{AWSC}) and the asymptotic neighbourhood condition (\ref{ANC}). Since the proof of the Hausdorff dimension relies on the Feng-Hu formula, it inherits its dependence on the structure of the dimensions of the principal-projection IFSs. To control these, we also require separation assumptions for the principal projections. See also the assumptions of \cite[Theorem 1.3]{rapaport2023dimensiondiagonalselfaffinesets}. 


\subsubsection{Hausdorff dimension}
Our first theorem asserts that under the assumption that the principal projections of the IFS satisfy the WSC, the Hausdorff dimension of \ref{A1} equals the limit of the maxima of the Bara\'nski formula applied to the functions at the $n$th level. 

For a finite set $A$, let $P(A):=\{p=(p_a)_{a\in A}\, \big|\, p_a\ge 0,\, \sum_{a\in A}p_a=1 \}$ be the set of probability vectors over $A$. Let $p \in P(\mathbb{G})$, $\ell \in \{1,2\}$. For $R\in\mathrm{p}_\ell\mathbb{G}$, let $q_{R}^\ell=\sum_{S\in\mathbb{G}: \mathrm{p}_\ell S=R}p_S$. For a \ref{A1} $\mathbb{G}$ and probability vector $p \in P(\mathbb{G})$, let

\begin{equation} \label{DpG}
	D(p,\mathbb{G}) := \begin{cases}
		\displaystyle
		\frac{\sum\limits_{S \in \mathbb{G}} p_S \log p_S}{\sum\limits_{S \in \mathbb{G}} p_S \log r_{S,1}}
		+ \frac{\sum\limits_{R \in \mathrm{p}_2\mathbb{G}} q_{R}^2 \log q_{R}^2}{\sum\limits_{S \in \mathbb{G}} p_S \log r_{S,2}}
		- \frac{\sum\limits_{R \in \mathrm{p}_2\mathbb{G}} q_{R}^2 \log q_{R}^2}{\sum\limits_{S \in \mathbb{G}} p_S \log r_{S,1}},  & \text{if} \quad \sum\limits_{S \in \mathbb{G}} p_S \log\frac{r_{S,2}}{ r_{S,1}} \geq 0 \\[25pt]
		\displaystyle
		\frac{\sum\limits_{S \in \mathbb{G}} p_S \log p_S}{\sum\limits_{S \in \mathbb{G}} p_S \log r_{S,2}}
		+ \frac{\sum\limits_{R \in \mathrm{p}_1\mathbb{G}} q_{R}^1 \log q_{R}^1}{\sum\limits_{S \in \mathbb{G}} p_S \log r_{S,1}}
		- \frac{\sum\limits_{R \in\mathrm{p}_1\mathbb{G}} q_{R}^1 \log q_{R}^1}{\sum\limits_{S \in \mathbb{G}} p_S \log r_{S,2}},  & \text{if} \quad \sum\limits_{S \in \mathbb{G}} p_S \log\frac{r_{S,2}}{r_{S,1}} < 0.
	\end{cases}
\end{equation}
Throughout the paper, we use the convention: $0\log 0:=0$. We define
$$
\mathbb{H}_{BA}(\mathbb{G})=\max\{D(p,\mathbb{G})\, \big|\, p \in P(\mathbb{G})\}.
$$

\begin{theorem}[Hausdorff dimension] \label{th}
Let $\mathbb{G}$ be an \ref{B2} \ref{A1} such that either the WSC holds for both projections, or it satisfies the \ref{G1}, the \ref{C1}, and the WSC for the $y$-projection. Then
\begin{equation}\label{eq:limdim}
\dH(\Lambda) = \lim_{n \to \infty} \mathbb{H}_{BA}(\mathbb{G}_n),
\end{equation}
where $\mathbb{G}_n$ is the set of $n$-fold compositions of the functions of $\mathbb{G}$.
\end{theorem}

Unfortunately, \eqref{eq:limdim} appears difficult to calculate for general weakly separated systems. For homogeneous systems, the formula simplifies.

\begin{corollary}[Hausdorff dimension for homogeneous systems] \label{T3}
	Let $\mathbb{G}$ be a \ref{B1}, \ref{B2} \ref{A1} satisfying the \ref{C1}, and the WSC for the $y$-projection with contraction ratios $0<r_1 < r_2<1$. Then
	\begin{equation*}
		\dH(\Lambda) = \lim_{n \to \infty}\frac{\log\left(\sum_{T \in \mathrm{p}_2\mathbb{G}_n} \#\{S\in\mathbb{G}_n\, \big|\, \mathrm{p}_{2} S=T\}^{\frac{\log r_2}{\log r_1}}\right)}{{-}n\log r_2}.
	\end{equation*}
\end{corollary}

\subsubsection{Box-counting dimension}
We assert that, under the assumption that the weak separation condition holds for the projections of the IFS, the box-counting dimension exists, and equals the limit of the maxima of the Feng-Wang formula applied to the functions at the $n$th level.

\begin{theorem}[Box-counting dimension] \label{tb}
Let $\mathbb{G}$ be a \ref{A1} such that either the WSC holds for both projections, or it satisfies the \ref{G1}, the \ref{C1}, and the WSC for the $y$-projection. For every $n\in\N$, let $d_n^1$ and $d_n^2$ be the unique real solutions to the equations
\begin{equation}\label{eq:boxdim}
1= \sum_{S\in \mathbb{G}_n} \left( \frac{|r_{S,1}|}{|r_{S,2}|} \right)^{\dB(\mathrm{p}_1\Lambda)} |r_{S,2}|^{d_n^1}, \qquad 1= \sum_{S \in \mathbb{G}_n} \left( \frac{|r_{S,2}|}{|r_{S,1}|} \right)^{\dB(\mathrm{p}_2\Lambda)} |r_{S,1}|^{d_n^2},
\end{equation}
Then
\begin{align*}
\dB(\Lambda) = \limsup_{n \to \infty} \max\{d_n^1,d_n^2\}.
\end{align*}
Furthermore, if for a fixed $\ell \in \{1,2\}, \ \mathrm{p}_\ell\mathbb{G}$ satisfies the WSC, then the definition of $d_n^\ell$ remains valid if \eqref{eq:boxdim} is replaced by 
\[
1= \sum_{S\in \mathbb{G}_n} \left( \frac{|r_{S,1}|}{|r_{S,2}|} \right)^{s_n^1} |r_{S,2}|^{d_n^1}, \qquad 1= \sum_{S \in \mathbb{G}_n} \left( \frac{|r_{S,2}|}{|r_{S,1}|} \right)^{s_n^2} |r_{S,1}|^{d_n^2},
\]
where $s_n^\ell$ is the unique real solution of the equation $\sum_{R\in\mathrm{p}_\ell\mathbb{G}_n}|r_{R,\ell}|^{s_n^\ell}=1$.
\end{theorem}

This theorem corresponds to a limiting version of the Feng-Wang formula \cite{MR2128947} for the box-counting dimension. Note that with respect to the box-counting dimension, the case of a \ref{A1} satisfying \ref{G1}, the \ref{C1}, and the WSC for the $y$-projection is a corollary of the results of Zerner \cite{MR1343732} and Feng and Wang \cite{MR2128947}, unlike the case of the Hausdorff dimension.

We establish explicit limit formulas that characterise the box-counting dimension. Finally, we state the corollary for homogeneous systems.
\begin{corollary}[Box-counting dimension for homogeneous systems] \label{tbh}
Let $\mathbb{G}$ be a \ref{B1} \ref{A1} satisfying the \ref{C1}, and the WSC for the $y$-projection. Then
\begin{equation*}
	\dB(\Lambda) = \lim_{n \to \infty}\frac{\log \big(\#\mathbb{G}_n\big)}{-n\log r_1} + \dB(\mathrm{p}_2\Lambda)\Big(1 - \frac{\log r_2}{\log r_1}\Big)
\end{equation*}
where $\dB(\mathrm{p}_2\Lambda)$ is the box-counting dimension of the attractor of the IFS $\mathrm{p}_2\mathbb{G}$.
\end{corollary}

The proofs differ in detail, but follow the same pattern: the upper bounds will be straightforward, while for the lower bound, we will need to find a close to \ref{B1} subsystem (an IFS whose functions are elements of $\mathbb{G}^*$), and then find a well-separated subsystem of that, with sufficiently large dimension. For the Hausdorff dimension, by Ferguson, Jordan, and Shmerkin \cite{ferguson2009hausdorffdimensionprojectionsselfaffine}, we can find a \ref{B1} subsystem for the first step; for the box-counting dimension, we will use a more elementary variant.

\subsection{Examples}

We now illustrate our results with two examples. The calculations of the following values are lengthy; therefore, we defer the proofs of these assertions to the final section of the paper.

\begin{example} \label{ex:1} Let $\mathbb{G}$ be the IFS depicted in Figure \ref{EF2}.

\begin{figure}[H]
	\centering
	\begin{minipage}{0.65\textwidth}
		\centering
		\scalebox{1}[1]{\includestandalone[width=0.75\textwidth]{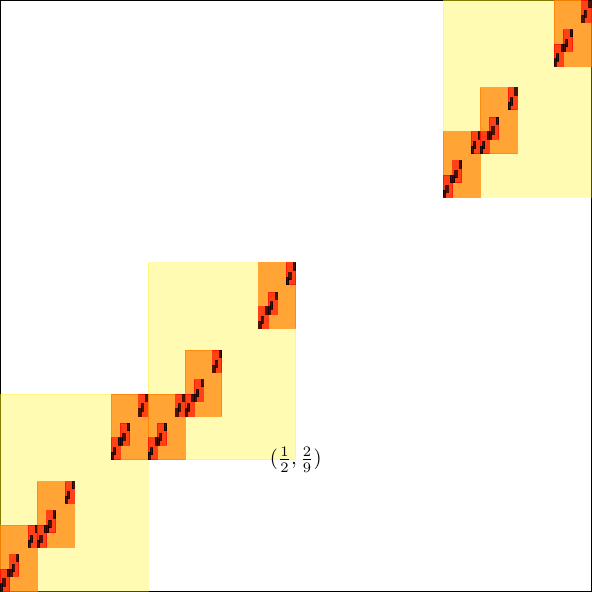}}
	\end{minipage}%
	\hfill
	\begin{minipage}{0.35\textwidth}
		\raggedright
		\begin{align*}
			\mathbb{G} := &\big\{S_1, S_2, S_3\big\}, \notag\\[6pt]
			S_1(x,y) &:=\left(\frac{x}{4},\frac{y}{3}\right), \notag\\[5pt]
			S_2(x,y) &:= \left(\frac{x}{4}+\frac14,\frac{x}{3}+\frac29\right), \notag\\[5pt]
			S_3(x,y) &:= \left(\frac{x}{4}+\frac34,\frac{y}{3}+\frac23\right). \notag
		\end{align*}
	\end{minipage}
	\caption{Example for the application of Corollary \ref{T3}.}
	\label{EF2}
\end{figure}
	
Then $\dH(\Lambda) = \log_3 \lambda^*$, where $\lambda^*$ is the unique $\lambda> 1$ such that
\begin{equation*}
\lambda = \frac{1}{\lambda-1}\sum_{k=2}^\infty k^\alpha \lambda^{2-k} + \frac{\lambda^2}{(\lambda-1)^3},
\end{equation*}
where $\alpha := \log_43$, and therefore
\begin{equation*}\begin{aligned}
\dH(\Lambda) &= \log_3 2.8960013515886529426596184724862681808317981559701975\dots \\
& = 0.967885533595539319438037445903385862252724017052009287837\dots
\end{aligned}\end{equation*}
\end{example}

\begin{example} \label{ex:2}
Let $\mathbb{G}$ be the IFS depicted in Figure \ref{EF1}.

\begin{figure}[H]
	\centering
	\begin{minipage}{0.65\textwidth}
		\centering
		\scalebox{1}[1]{\includestandalone[width=0.75\textwidth]{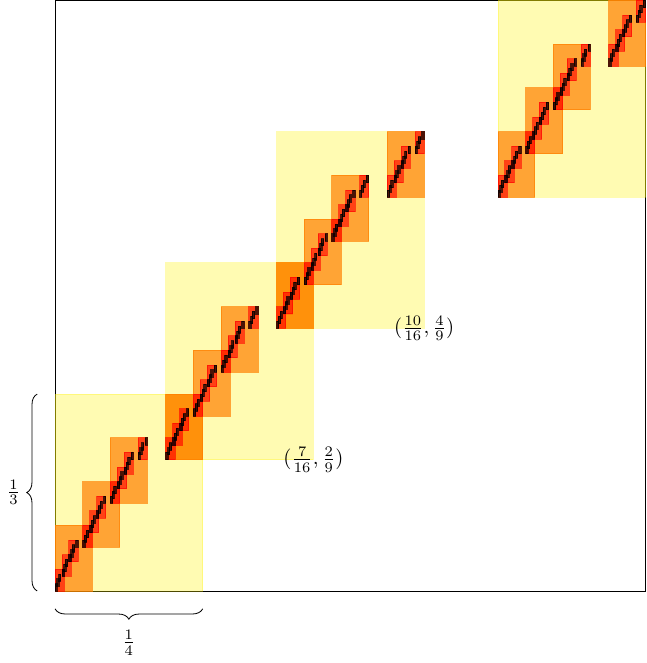}}
	\end{minipage}%
	\hfill
	\begin{minipage}{0.35\textwidth}
		\raggedright
		\begin{align*}
			\mathbb{G} := &\big\{S_1, S_2, S_3, S_4 \big\}, \notag\\[6pt]
			S_1(x,y) &:= \left(\frac{x}{4},\frac{y}{3}\right), \notag\\[5pt]
			S_2(x,y) &:= \left(\frac{x}{4}+\frac{3}{16},\frac{y}{3}+\frac29\right), \notag\\[5pt]
			S_3(x,y) &:= \left(\frac{x}{4}+\frac{6}{16},\frac{y}{3}+\frac49\right), \notag\\[5pt]
			S_4(x,y) &:= \left(\frac{x}{4}+\frac34,\frac{y}{3}+\frac23\right). \notag
		\end{align*}
	\end{minipage}
	\caption{Example for the application of Corollary \ref{tbh}.}
	\label{EF1}
\end{figure}

Then
\begin{align*}
\dB(\Lambda) = \frac{\log 3}{\log 3}\Big(1 - \frac{\log 3}{\log 4}\Big) + \frac{\log(2+\sqrt{2})}{\log 4} = 1.093295401221 \dots
\end{align*}
while $\dH(\Lambda) = \log_3 \lambda^*$, where $\lambda^*$ is the unique $\lambda> 1$ such that
\begin{equation*}
\lambda^2 = 2\lambda + (\lambda-1)\sum_{k=0}^\infty (k+1)^\alpha \lambda^{-k},
\end{equation*}
where $\alpha := \log_43$, and therefore
\begin{align*}
\dH(\Lambda) &= \log_3 3.3053444391403030804351198628168101723048873990180\dots \\
& = 1.08822802893857463333338129344577152156663085706579\dots
\end{align*}

\end{example}

We note that the examples above satisfy the assumptions of He, Lau, and Rao \cite[Theorem~4.4]{He2003} based on Kenyon and Peres \cite[Theorem~2.2]{Kenyon1996}. However, we provide an alternative method for calculating the dimension values.

\subsubsection*{Acknowledgement}

We wish to express our gratitude to the anonymous referee for their leniency and their valuable comments, especially for pointing out the result of Ferguson, Jordan, and Shmerkin \cite{ferguson2009hausdorffdimensionprojectionsselfaffine}, which shortened the previous argument for the Hausdorff dimension significantly.

\section{Preliminaries}

\subsection{Notations}
\label{subsection:symbolic}
Let $\mathbb{G}=\{S_i\}_{ i \in \Sigma}$ be an IFS indexed by $\Sigma = \{1,\dots,\#\mathbb{G}\}$. For simplicity, we also write for the contraction ratios $r_{i,\ell}:= r_{S_i,\ell}$. For $k\in \N^+$, let $\Sigma^k$ be the set of $k$-tuples formed by the elements of $\Sigma$. Let $\Sigma^*$ be the set $\bigcup_{k \in \N^+}\Sigma^k$. Finite words in $\Sigma^*$ will be denoted in mathfrak style: $\iii,\jjj, \kkk, \hhh$, except for words of length 1, in which case we might use $i,j,k,h$.
	
For a finite word $\iii\in\Sigma^*$, denote by $|\iii|$ the length of $\iii$. If $k<\ell \in \N$ and $\iii$ is a finite word with $\ell\leq|\iii|$, we denote:
	\begin{align*}
		\iii|_{(k,\ell]} &:= (\iii_{k+1} \dots \iii_\ell) \\
		\iii- &:=  (\iii_{1} \dots \iii_{|\iii|-1}) = \iii|_{(0,|\iii|-1]}.
	\end{align*}

For $\iii=(i_1,\ldots,i_{|\iii|}) \in \Sigma^*$, write
\begin{align*}
	S_{\iii} := S_{i_1} \circ \dots \circ S_{i_{|\iii|}}.
\end{align*}

The above symbolic space describes separated sets very well, but when less is assumed for separation, it shows some shortcomings. Namely, if there are exact overlaps, many finite words are redundant.

For $n \in \N^+$ let the $n$th level functions be
\begin{align*}
\mathbb{G}_n := \big\{ S_{\iii} \, \big|\, \iii \in \Sigma^n \big\}.
\end{align*}
We define $\Sigma^{\{n\}} \subseteq \Sigma^n$ iteratively. Let $\Sigma^{\{1\}}:=\Sigma$. Given $\Sigma^{\{k\}}$, we define $\Sigma^{\{k+1\}}$ as a maximal subset of $\Sigma^{k+1}$ such that
\begin{align*}
\mathbb{G}_{k+1}=\Big\{S_{\iii}\, \Big|\, \iii\in\Sigma^{\{k+1\}}\Big\}, S_{\iii}\neq S_{\jjj}\text{ for $\iii\neq\jjj\in\Sigma^{\{k+1\}}$ and }\Sigma^{\{k\}} \supseteq \Big\{ \iii- \, \Big|\, \iii \in \Sigma^{\{k+1\}} \Big\}
\end{align*}
for every $k\in\N$.
By construction, there is a one-to-one correspondence between $\Sigma^{\{n\}}$ and the maps in $\mathbb{G}_n$. Similarly, for $\ell\in\{1,2\}$, let $\Gamma_\ell^{\{n\}}\subseteq\Sigma^{\{n\}}$ be such that
\begin{align*}
\mathrm{p}_\ell\mathbb{G}_{n}=\Big\{\mathrm{p}_\ell S_{\iii}\, \Big|\, \iii\in\Gamma_\ell^{\{n\}}\Big\}, \mathrm{p}_\ell S_{\iii}\neq \mathrm{p}_\ell S_{\jjj}\text{ for $\iii\neq\jjj\in\Gamma_\ell^{\{n\}}$} \text{ and }\Gamma_\ell^{\{n-1\}} \supseteq \Big\{ \iii- \, \Big|\, \iii \in \Gamma_\ell^{\{n\}} \Big\}.
\end{align*}
There is a one-to-one correspondence between the elements of $\Gamma_\ell^{\{n\}}$ and the maps in $\mathrm{p}_\ell\mathbb{G}_n$. Let $\Gamma_\ell := \Gamma_\ell^{\{1\}}$. For $\iii \in \Sigma^n$, let $\mathrm{p}_\ell (\iii)$ be the unique element of $\Gamma_\ell^{\{n\}}$ such that
\begin{align*}
\mathrm{p}_\ell S_{\mathrm{p}_\ell (\iii)} = \mathrm{p}_\ell S_{\iii}.
\end{align*}
Let $\Sigma^{\{*\}}=\bigcup_{n=1}^\infty\Sigma^{\{n\}}$. Furthermore, for $m \in \N^+$, let
\begin{align*}
\Sigma^{\{n\}m} &:= \Big\{ \iii^1\iii^2 \dots \iii^m \, \Big|\, \iii^j \in \Sigma^{\{n\}},\ j \in \{1,\dots,m\} \Big\}, \quad \Sigma^m := \Sigma^{\{1\}m}, \\
\Gamma_\ell^{\{n\}m} &:=\Big\{\jjj^1\jjj^2\dots \jjj^m \, \Big| \, \jjj^k\in\Gamma_\ell^{\{n\}},\ k \in \{1,\dots,m\}\Big\}, \quad \Gamma_\ell^m := \Gamma_\ell^{\{1\}m}.
\end{align*}

We now introduce some notation specific to diagonally aligned self-affine sets. For $\ell\in\{1,2\}$, let
\begin{align*}
r_{\max,\ell} &:= \max_{i \in \Sigma} \big\{ |r_{i,\ell}| \big\} &  r_{\min,\ell} := &\min_{i \in \Sigma} \big\{ |r_{i,\ell}| \big\}\\
r_{\max} &:= \max\big\{ r_{\max,1},r_{\max,2} \big\} & r_{\min} := &\min\big\{r_{\min,1},r_{\min,2} \big\}.
\end{align*}
For
$\iii = (i_1,i_2,\dots, i_{|\iii|}) \in \Sigma^*$, write
\[ r_{\iii,1}:= \prod_{m=1}^{|\iii|} r_{i_{m},1}\quad \quad r_{\iii,2}:= \prod_{m=1}^{|\iii|} r_{i_{m},2}.\]

Define the $\delta$-Moran cut-set of the self-similar IFS $\mathbb{G} = \big\{S_i \big\}_{i \in \Sigma}$ and $\delta \in (0,1)$ as
\[
M_\delta(\mathbb{G}, \Sigma) := \Big\{ \iii\in\Sigma^{\{*\}} \, \Big|\, |S_{\iii}(\Lambda)| \leq \delta < |S_{\iii-}(\Lambda)| \Big\}.
\]
Furthermore, define
\[
\mathcal{M}_\delta(\mathbb{G}, \Sigma) := \Big\{ S_{\iii}\, \Big|\, \iii \in M_\delta(\mathbb{G}, \Sigma)\Big\}.
\]
In particular, given a diagonally aligned self-affine set $\mathbb{G}_n$, for $\ell \in \{1,2\}$ denote $M_\delta(\mathrm{p}_\ell\mathbb{G}, \Gamma_\ell)$ by $M_\delta^\ell$, and $\mathcal{M}_\delta(\mathrm{p}_\ell\mathbb{G}, \Gamma_\ell)$ by $\mathcal{M}_\delta^{\ell}$. Let
\begin{align} \label{eq:Delta}
	\Delta_\delta := \Big\{\iii \in \Sigma^{\{*\}} \, \Big|\, \min\{|r_{\iii,1}|,|r_{\iii,2}|\} \leq \delta < \min\{|r_{\iii-,1}|,|r_{\iii-,2}|\}\Big\} \quad \text{ and } \quad \mathbb{S}_\delta := \big\{S_{\iii}  \, \big|\, \iii \in \Delta_\delta \big\}.
\end{align}
We define
\begin{align*}
\Delta_\delta^1 := \big\{\iii \in \Delta_\delta \, \big|\, |r_{\iii,1}| > |r_{\iii,2}|\big\}, \quad \mathbb{S}_\delta^1 := \big\{S_{\iii}  \, \big|\, \iii \in \Delta_\delta^1 \big\},\\
	\Delta_\delta^2 := \big\{\iii \in \Delta_\delta \, \big|\, |r_{\iii,2}| \geq |r_{\iii,1}|\big\}, \quad \mathbb{S}_\delta^2 := \big\{S_{\iii}  \, \big|\, \iii \in \Delta_\delta^2 \big\} .
\end{align*}

\subsection{Separation conditions}
Given a self-similar IFS $\mathbb{G} = \{S_i\}_{i \in \Sigma}$, let
\begin{align} \label{eq:awscdefqu}
g(\delta) := \sup_{\gamma \in [\delta,1]}\sup_{x \in \R^2}\#\left\{S \in \mathcal{M}_{\gamma}(\mathbb{G}, \Sigma) \, \middle|\, S(\Lambda) \cap B(x,\gamma) \neq \emptyset \right\}
\end{align}

\begin{enumerate}[label = $$, leftmargin=0pt, ref = AWSC]
\item\label{AWSC} We say that a self-similar IFS $\mathbb{G}$ satisfies the \textbf{asymptotic weak separation condition} (AWSC) if $g(\delta)$ is subpolynomial, i.e.
\begin{align*}
\lim_{\delta \to 0}\frac{ \log g(\delta)}{-\log \delta} = 0.
\end{align*}
\end{enumerate}

Feng \cite{MR2322179} defined the AWSC in a slightly different form, without the $\sup_{\gamma \in [\delta,1]}$ part. However, that only ensures the monotonicity of $g$. We remark that although for a self-similar IFS on the line, the WSC is atypical, the \ref{AWSC} is typical in some sense. Given a self-similar IFS satisfying the exponential separation condition, the \ref{AWSC} holds if and only if the similarity dimension is $\leq 1$. See Barral and Feng \cite{barral2020multifractalformalismselfsimilarmeasures}.

In Anttila, Bárány, and Käenmäki \cite{ANTTILA_2023}, the authors define a separation condition for self-affine sets, called the bounded neighbourhood condition. This holds if there exists a constant $C$ such that any ball of radius $\delta$ is intersected by at most $C$ cylinders from the functions with the smallest contraction ratio approximately $\delta$.

Similarly, we define its asymptotic version for \ref{A1}s.

\begin{enumerate}[label = $$, leftmargin=0pt, ref = ANC]
\item\label{ANC} Let $\mathbb{G}$ be a \ref{A1}. We say that the \textbf{asymptotic neighbour condition} (ANC) holds if
\[
f(\delta) := \sup_{\gamma \in [\delta,1]}\sup_{x \in \R^2}\#\left\{S \in \mathbb{S}_\gamma \, \middle|\, S(\Lambda) \cap B(x,\gamma) \neq \emptyset \right\}
\]
is subpolynomial, i.e. $\log f(\delta)/\log\delta\to0$ as $\delta\to0$.
\end{enumerate}

We note that, for a \ref{A1}, satisfying the \ref{ANC} does not imply the \ref{AWSC} for its principal projections, but the reversed implication hold. 
\begin{lemma}
Let $\mathbb{G}$ be a \ref{A1} such that its two principal projections satisfy the \ref{AWSC}, then $\mathbb{G}$ satisfies the \ref{ANC}.
\end{lemma}
\begin{Proof}

Let $x \in \R^2, \delta > 0$. By symmetry, it is sufficient to bound the asymptotics of the cardinality of 
\begin{align*}
H(x,\delta) := \left\{S \in \mathbb{S}_\delta^1 \, \middle|\, S(\Lambda) \cap B(x,\delta) \neq \emptyset \right\}.
\end{align*}
Since $\delta r_{\min} < |r_{S,1}|\leq \delta^{\frac{\log r_{\max}}{\log r_{\min}}}$ for any $S\in\mathbb{S}_\delta^1$,
\begin{multline*}
\#H(x,\delta) \leq \#\left\{\mathrm{p}_1S \in \mathcal{M}_\delta^1 \, \middle|\, (\mathrm{p}_1S)(\Lambda) \cap B(\mathrm{p}_1(x),\delta) \neq \emptyset \right\} \\ \cdot \#\left\{\mathrm{p}_2S \, \middle|\, \delta r_{\min} < |r_{S,1}|\leq \delta^{\frac{\log r_{\max}}{\log r_{\min}}} , \, (\mathrm{p}_2S)(\Lambda) \cap B(\mathrm{p}_2(x),\delta) \neq \emptyset\right\}.
\end{multline*}
The first term grows subpolynomially since $\mathrm{p}_1\mathbb{G}$ satisfies \ref{AWSC}; thus, we need to bound the second term. Notice that
\begin{multline*}
\left\{\mathrm{p}_2S \, \middle|\, \delta r_{\min} < |r_{S,1}|\leq \delta^{\frac{\log r_{\max}}{\log r_{\min}}} , \, (\mathrm{p}_2S)(\Lambda) \cap B(\mathrm{p}_2(x),\delta) \neq \emptyset\right\} \\ \subseteq \bigcup_{n=0}^h \left\{\mathrm{p}_2S \in \mathcal{M}_{\delta r_{\max}^{-n}}^2 \, \middle|\, (\mathrm{p}_2S)(\Lambda) \cap B(\mathrm{p}_2(x),\delta) \neq \emptyset\right\},
\end{multline*}
where 
\begin{align*}
h := \left\lceil \frac{\log \delta}{\log r_{\max}}\Big(1-\frac{\log r_{\max}}{\log r_{\min}}\Big)\right\rceil.
\end{align*}
Since $\mathrm{p}_2\mathbb{G}$ satisfies the \ref{AWSC}, there is a subpolynomial function $g(\delta)$ such that for any $n \in \N$ and any $y \in \R$
\begin{align*}
\#\left\{\mathrm{p}_2S \in \mathcal{M}_{\delta r_{\max}^{-n}}^2 \middle|\, (\mathrm{p}_2S)(\Lambda) \cap B(y,\delta) \neq \emptyset\right\} \leq g(\delta r_{\max}^{-n}) \leq g(\delta).
\end{align*}
Thus
\begin{align*}
\# \left(\bigcup_{n=0}^h \left\{\mathrm{p}_2S \in \mathcal{M}_{\delta r_{\max}^{-n}}^2 \middle|\, (\mathrm{p}_2S)(\Lambda) \cap B(\mathrm{p}_2(x),\delta) \neq \emptyset\right\} \right) \leq h \cdot g(\delta) =: g_1(\delta),
\end{align*}
where $g_1(\delta)$ is a subpolynomial function.
\end{Proof}

\subsection{Dimension estimates for carpets}

We now state some well-known results from the dimension theory of \ref{A1}s.

\begin{theorem}[Feng and Wang \cite{MR2128947}]\label{thm:FW}
	Let $\mathbb{G}$ be a \ref{A1} satisfying the \ref{C1}. Then \[\dB(\Lambda)=\max\{d_1,d_2\},\]where $d_1$ and $d_2$ are the unique reals defined by the equations
	\[
	1= \sum_{S\in \mathbb{G}} \left( \frac{|r_{S,1}|}{|r_{S,2}|} \right)^{\dB(\mathrm{p}_1\Lambda)} |r_{S,2}|^{d_1} \quad \text{and} \quad 1= \sum_{S \in \mathbb{G}} \left( \frac{|r_{S,2}|}{|r_{S,1}|} \right)^{^{\dB(\mathrm{p}_2\Lambda)}} |r_{S,1}|^{d_2}.
	\]
\end{theorem}

Let us note that Theorem \ref{thm:FW} is not stated in this way in Feng and Wang \cite{MR2128947}. For a proof of this formula, we refer to \cite[Theorem~11.4.2]{MR4661364}.

\begin{theorem}[Lalley and Gatzouras \cite{MR1183358}]\label{thm:GL}
	Let $\mathbb{G}$ be a \ref{A1} satisfying \ref{B2}, the \ref{C1}, \ref{G1} and the OSC for the $y$-projection. Then
	\begin{equation*}
		\dH(\Lambda) = 	\mathbb{H}_{BA}(\mathbb{G})=\max_{p\in P(\mathbb{G})}\left\{\frac{\sum\limits_{S \in \mathbb{G}} p_S \log p_S}{\sum\limits_{S \in \mathbb{G}} p_S \log r_{S,1}}
		+ \frac{\sum\limits_{R \in \mathrm{p}_2\mathbb{G}} q_{R}^2 \log q_{R}^2}{\sum\limits_{S \in \mathbb{G}} p_S \log r_{S,2}}
		- \frac{\sum\limits_{R \in \mathrm{p}_2\mathbb{G}} q_{R}^2 \log q_{R}^2}{\sum\limits_{S \in \mathbb{G}} p_S \log r_{S,1}}\right\},
	\end{equation*}
	where $q_{R}^2=\sum_{S\in\mathbb{G}\ :\ \mathrm{p}_2S=R}p_S$ and .
\end{theorem}

The proof can be found in \cite[Theorem~5.3]{MR1183358}.

\begin{theorem}[Pardo-Sim\'on \cite{simon2017dimensionsoverlappinggeneralizationbaranski}]\label{thm:PS}
	Let $\mathbb{G}$ be a \ref{A1}. Then
	\begin{equation*}
		\dH(\Lambda) \leq	\mathbb{H}_{BA}(\mathbb{G})=\max\{D(p,\mathbb{G})\, \big|\, p \in P(\mathbb{G})\}.
	\end{equation*}
\end{theorem}

For a proof, we refer to \cite[Proposition~3.5]{simon2017dimensionsoverlappinggeneralizationbaranski}.

\vspace{15 pt}

\section{Hausdorff dimension}

We now state a more general version of our result on the Hausdorff dimension, which clearly implies Theorem \ref{th}.

\begin{theorem}[Hausdorff dimension] \label{T2}
Let $\mathbb{G}$ be an \ref{B2} \ref{A1} such that either the \ref{AWSC} holds for both projections, or it satisfies the \ref{G1}, the \ref{ANC} and the \ref{AWSC} for the $y$-projection. Then
\begin{equation*}
\dH(\Lambda) = \lim_{n \to \infty} \mathbb{H}_{BA}(\mathbb{G}_n),
\end{equation*}
where $\mathbb{G}_n$ is the set of $n$-fold compositions of $\mathbb{G}$. 
\end{theorem}

\begin{remark}
    Let us note that from the proof of Theorem~\ref{T2}, a construction of a \ref{B1} subsystem, which satisfies separation conditions on the plane and in the principal direction, In particular, for any $\eps >0$ there exists $N \in \N$ such that for all $n>N$ there is $G \subseteq \mathbb{G}_n$, a \ref{B1}, \ref{B2} \ref{A1} satisfying the \ref{C1}, the OSC for one of the principal projections and with
\begin{equation*}
\dH(\Lambda) \leq \dH(\Lambda(G)) + \eps.
\end{equation*}
\end{remark}



\subsection{Homogeneity lemmas}

The AWSC asserts that the number of cylinders in a $\delta$-Moran cut-set that intersect a $\delta$-ball is subpolynomial in $\delta$. However, the projections of level-$n$ functions on a general self-affine carpet are not Moran cut-sets. Thus, we first want to relate our system to one with homogeneous contractions. This can be done with an $\eps > 0$ loss by a lemma of Ferguson, Jordan, and Shmerkin \cite{ferguson2009hausdorffdimensionprojectionsselfaffine}.

The statement of \cite[Lemma 4.3]{ferguson2009hausdorffdimensionprojectionsselfaffine} considers Bara\'nski and Lalley-Gatzouras systems, but their proof is more general. To find a sufficient subsystem $G$, they consider $\mathbb{G}_m$ for $m$ sufficiently large, and they select functions that attain the $i$th function of $\mathbb{G}$ over their $m$ components approximately $p_i \cdot m$ times in $G$. They show that then $\max_pD(p,G)$ is $\eps$ close to $\max_pD(p,\mathbb{G})$. To translate this to the overlapping case, we need to work with $\Sigma^m$ instead of $\Sigma^{\{m\}}$ and with $\Gamma^m_\ell$ instead of $\Gamma^{\{m\}}_\ell$ for $\ell \in \{1,2\}$. For the definitions of these notions, we direct the reader to Subsection \ref{subsection:symbolic}. 

Write $P:=P(\Sigma), P^n:=P(\Sigma^n), P^{\{n\}}:=P(\Sigma^{\{n\}})$.
We denote probability vectors by $p,q,s$. $p$ will be used as a probability distribution on the symbolic space of \ref{A1}, and $q$ and $s$ will be probability distributions defined on the symbolic space of the principal projections.

For $p \in P^m, \ell \in \{1,2\}$, let $s^\ell$ be distributions on $\Gamma_\ell^{m}$ defined as 
\[ 
s^{\ell}_{\jjj} := \sum_{\kkk \in \Sigma^m:\ \forall v = 1,\dots,m\ \mathrm{p}_\ell (\kkk_v)= \mathrm{p}_\ell (\jjj_v)}p_{\kkk}, \qquad \text{for  } \jjj \in \Gamma_\ell^{m}. 
\]
For $p \in P^m$, define 
\begin{equation*}\label{eq:defD}
	D'(p) := \begin{cases}
		\displaystyle
		\frac{\sum\limits_{\iii \in \Sigma^m} p_{\iii} \log p_{\iii}}{\sum\limits_{\iii \in \Sigma^m} p_{\iii} \log r_{\iii,1}}
		+ \frac{\sum\limits_{\jjj \in \Gamma_2^{m}} s^{2}_{\jjj} \log s^{2}_{\jjj}}{\sum\limits_{\iii \in \Sigma^m} p_{\iii} \log r_{{\iii},2}}
		- \frac{\sum\limits_{{\jjj} \in \Gamma_2^{m}} s^{2}_{\jjj} \log s^{2}_{\jjj}}{\sum\limits_{\iii \in \Sigma^m} p_{\iii} \log r_{{\iii},1}},  & \text{if} \quad \sum\limits_{\iii \in \Sigma^m} p_{\iii} \log\frac{r_{{\iii},2}}{ r_{{\iii},1}} \geq 0 \\[25pt]
		\displaystyle
		\frac{\sum\limits_{\iii \in \Sigma^m} p_{\iii} \log p_{\iii}}{\sum\limits_{\iii \in \Sigma^m} p_{\iii} \log r_{\iii,2}}
		+ \frac{\sum\limits_{{\jjj} \in \Gamma_1^{m}} s^{1}_{\jjj} \log s^{1}_{\jjj}}{\sum\limits_{\iii \in \Sigma^m} p_{\iii} \log r_{{\iii},1}}
		- \frac{\sum\limits_{{\jjj} \in \Gamma_1^{m}} s^{1}_{\jjj} \log s^{1}_{\jjj}}{\sum\limits_{\iii \in \Sigma^m} p_{\iii} \log r_{{\iii},2}},  & \text{if} \quad \sum\limits_{\iii \in \Sigma^m} p_{\iii} \log\frac{r_{{\iii},2}}{r_{{\iii},1}} < 0.
	\end{cases}
\end{equation*}

We remark that $D(p, \mathbb{G}) = D'(p)$ for $p \in P = P^1$. Furthermore, if $\mathrm{p}_\ell\mathbb{G}$ satisfies the OSC, then $D'(p) = D(p,\mathbb{G})$.

\begin{lemma}[Ferguson, Jordan, Shmerkin \cite{ferguson2009hausdorffdimensionprojectionsselfaffine}, Lemma 4.3] \label{lemma:hhomogen2}
Let $\mathbb{G} = \{S_i\}_{i \in \Sigma}$ be a \ref{A1}. Then for any $\varepsilon > 0$ there is $M \in \N$ such that for any $m > M$ there is a \ref{B1} \ref{A1}, $G = \{S_{\iii}\}_{\iii \in \Sigma'}$ such that $\Sigma' \subseteq \Sigma^m$ and 
\[ \max_{p \in P}D(p,\mathbb{G}) \leq  \max_{p \in P(\Sigma')}D'(p) + \varepsilon .\]
\end{lemma}

By the following lemma, we can compute the maximum of the Feng-Hu formula for a \ref{B1} \ref{A1}.

\begin{lemma}\label{lem:lagrange}
	Let $\mathbb{G} = \{S_i\}_{i \in \Sigma}$ be a \ref{B1} \ref{A1}. Then 
	\begin{multline*}
\max_{p\in P(\Sigma)}D(p, \mathbb{G}) =
	\max_{p\in P(\Sigma)}\left\{\frac{\sum\limits_{i \in \Sigma} p_i \log p_i}{\log r_1}
	+ \frac{\sum\limits_{j \in \Gamma_2} q_j^2 \log q_j^2}{\log r_2}
	- \frac{\sum\limits_{j \in \Gamma_2} q_j^2 \log q_j^2}{\log r_1}\right\} \\
=\frac{\log \Big(\sum\limits_{j\in\Gamma_2}\#\{i\in\Sigma\, |\, \mathrm{p}_2(i)=j\}^{\frac{\log r_2}{\log r_1}}\Big)}{-\log r_2}.
\end{multline*}
In particular, for the IFS $G = \{S_{\iii}\}_{\iii \in \Sigma'},\Sigma' \in \Sigma^m$ from Lemma \ref{lemma:hhomogen2}
\begin{multline*}
	\max_{p \in P(\Sigma')}D'(p) = \max_{p\in P(\Sigma')}\left\{\frac{\sum\limits_{{\iii} \in \Sigma'} p_{\iii} \log p_{\iii}}{\log r_1}
	+ \frac{\sum\limits_{{\jjj} \in \Gamma_2'} s_{\jjj}^2 \log s_{\jjj}^2}{\log r_2}
	- \frac{\sum\limits_{{\jjj} \in \Gamma_2'} s_{\jjj}^2 \log s_{\jjj}^2}{\log r_1}\right\} \\
=\frac{\log \Big(\sum\limits_{{\jjj}\in\Gamma_2'}\#\big\{{\iii}\in\Sigma'\, \big|\, \mathrm{p}_2({\iii}_k)={\jjj}_k,\, k = 1,\dots,m\big\}^{\frac{\log r_2}{\log r_1}}\Big)}{-\log r_2}.
\end{multline*}
where $\Gamma_2' \subseteq \Gamma_2^m$ and where $r_1 = r_{\iii,1}, r_2 = r_{\iii,2}$ for any $\iii \in \Sigma'$.
\end{lemma}

The proof follows by \cite[Proposition~3.4]{MR1183358}, or by using the Lagrange multiplier method.

\vspace{5pt}

\emph{Summary of the proof of Theorem~\ref{T2}:}
To prove Theorem \ref{T2}, we first consider the $n$-fold functions, $\mathbb{G}_n$ (and its symbolic space $\Sigma^{\{n\}}$) for a fixed $n \in \N$. Then, by Lemma~\ref{lemma:hhomogen2}, for any $\varepsilon > 0$ there exists $M>0$ such that for any $m > M$ with respect to the Feng-Hu formula, $\Sigma^{\{n\}}$ will be $\eps$ close to a homogeneous subsystem of $\Sigma^{\{n\}m}$. Now, using the asymptotic weak separation, we will see that this introduces exponentially many overlaps in $m$, and a subexponential number of overlaps in $n$. Finally, we will take $m \to \infty$ so that $\eps \to 0$, and then $n \to \infty$ so that the constant arising from the exponential number of overlaps would tend to $0$.

\subsection{Lemmas working with the asymptotic weak separation condition}
Here, we derive a covering lemma from our separation assumptions. Let $\mathbb{G}$ be a \ref{A1}. Recall $\Sigma^{\{n\}m}$ and $\Gamma_\ell^{\{n\}m}$ for $n,m \in \N^+, \ell \in \{1,2\}$.

\begin{lemma}\label{lemma:symbolichomogene}
For any $\eps > 0$ and $n \in \N^+$, there is $M \in \N$ such that for any $m>M$ there exists $\Sigma' \subseteq \Sigma^{\{n\}m}$ and $\ell \in \{1,2\}$ such that
the IFS $\{S_{\iii}\}_{\iii \in \Sigma'}$ is a \ref{B1} \ref{A1} with common contraction ratios $r_1',r_2' \in (0,1)$ with $r_{\jjj,\ell} = r_\ell' \geq r_{3-\ell}' = r_{\jjj,3-\ell}$ for any $\jjj\in\Sigma'$; moreover,
\begin{align*}
\max_{p \in P^{\{n\}}}D(p,\mathbb{G}_n) \leq \frac{\log \bigg(\sum\limits_{\jjj\in\Gamma_{\ell}^{\{n\}m}}(\#\Sigma_{\jjj}^\ell )^{\frac{\log r_\ell'}{\log r_{3-\ell}'}}\bigg)}{-\log r_\ell'} + \eps, 
\end{align*}
where,  
\begin{align} \label{eq:projrow}
\Sigma_{\jjj}^\ell := \big\{\iii \in \Sigma' \, \big|\, \ \mathrm{p}_{\ell} (\iii |_{(m(v-1),mv]})= \mathrm{p}_{\ell} (\jjj |_{(m(v-1),mv]})\,\, \forall v = 1,\dots,m\big\}\quad\text{ for $\jjj \in \Gamma_\ell^{\{n\}m}$}.
\end{align}
\end{lemma}

\begin{Proof}
The claim follows by applying Lemma \ref{lemma:hhomogen2} to the IFS $\mathbb{G}_n$, and then by applying Lemma \ref{lem:lagrange}.
\end{Proof}

For simplicity, let us introduce a notation for the restricted projection $\mathrm{p}_\ell'\colon\Sigma'\to\Gamma_\ell^{\{n\}m}$ by
\begin{align*}
\mathrm{p}'_\ell(\iii) = \mathrm{p}'_\ell(\jjj) \quad \text{ if and only} \quad \forall v = 1,\dots,m: \ \mathrm{p}_{\ell} (\iii |_{(m(v-1),mv]})= \mathrm{p}_{\ell} (\jjj |_{(m(v-1),mv]}).
\end{align*}

Let $\Gamma_\ell' \subseteq \Sigma'$ be minimal such that $\{\mathrm{p}_\ell'(\iii) \, | \, \iii \in \Sigma'\} = \{\mathrm{p}_\ell'(\iii) \, | \, \iii \in \Gamma_\ell'\}.$ We note that $\Gamma_\ell'$ might not be a subset of $\Gamma_\ell^{\{n\}m}$.

We say that a function $C(n)$ is subexponential if 
\begin{align*}
\lim_{n \to \infty}\frac{1}{n} \log C(n) = 0.
\end{align*}

\begin{lemma} \label{lemma:WSCnew}
There exists subexponential function $C(n)$ such that for any $n,m \in \N^+$ and for any $\Sigma' \subseteq \Sigma^{\{n\}m}$ such that $\{S_{\iii}\}_{\iii \in \Sigma'}$ is a \ref{B1} \ref{A1},
\begin{itemize}
\item if $\mathrm{p}_\ell \mathbb{G}$ satisfies the \ref{AWSC}, then 
\begin{align}
\forall \jjj \in \Gamma_{\ell}': \quad \# \big\{\iii \in \Gamma_\ell' \, \big|\, (\mathrm{p}_\ell S_{\iii})((0,1)) \cap (\mathrm{p}_\ell S_{\jjj})((0,1)) \neq \emptyset \big\} &\leq C(n)^{m},\label{eq:WSCuse2}
\end{align}
\item if $\mathbb{G}$ satisfies the \ref{ANC} then
\begin{align}
\forall \jjj \in \Sigma': \quad \# \left\{\iii \in \Sigma_{\jjj}' \, \middle|\, S_{\iii}((0,1)^2) \cap S_{\jjj}((0,1)^2) \neq \emptyset \right\} &\leq C(n)^{m}\label{eq:WCSuse1}.
\end{align}
\end{itemize}
\end{lemma}

\begin{Proof}

Let us first consider the first statement. Clearly,
\begin{align*}
\Gamma_\ell^{\{n\}} \subseteq \bigcup_{k = 0}^K M_{r_{\max}^{n+k}}^\ell, \quad \text{where} \quad K := \left\lceil n \Big( \frac{\log r_{\min}}{\log r_{\max}}-1\Big) \right\rceil.
\end{align*}
Let $\mathrm{p}_\ell \mathbb{G}$ satisfy the \ref{AWSC} with bounding function $g(\delta)$. Then for any $x \in \R$,
\begin{align*}
\# \big\{\iii \in \Gamma_\ell' \, \big|\, x \in (\mathrm{p}_\ell S_{\iii})((0,1)) \big\} &\leq \# \big\{\iii \in \Gamma_\ell^{\{n\}m} \, \big|\, x \in (\mathrm{p}_\ell S_{\iii})((0,1))\big\} \\ &\leq \left(\sup_{x\in[0,1]}\# \big\{\iii \in \Gamma_\ell^{\{n\}} \, \big|\, x \in (\mathrm{p}_\ell S_{\iii})((0,1)) \big\}\right)^m \leq (K \cdot g(r_{\min}^n))^m,
\end{align*}
where, for the last inequality, we used the definition of \ref{AWSC} for balls centred at $x$ in each Moran cut-set.
Now, using the homogeneity of $\Gamma_{\ell}'$, i.e. that $|(\mathrm{p}_\ell S_{\iii})((0,1))| = |(\mathrm{p}_\ell S_{\jjj})((0,1))|$ for any $\iii, \jjj \in \Gamma_{\ell}'$, if $(\mathrm{p}_\ell S_{\iii})((0,1)) \cap (\mathrm{p}_\ell S_{\jjj})((0,1)) \neq \emptyset$, then $(\mathrm{p}_\ell S_{\iii})((0,1))\cap\{(\mathrm{p}_\ell S_{\jjj})(0),(\mathrm{p}_\ell S_{\jjj})(1/2),(\mathrm{p}_\ell S_{\jjj})(1)\}\neq\emptyset$. Thus, for any $\jjj\in\Gamma_\ell'$,
\begin{align*}
\# \big\{\iii \in \Gamma_\ell' \, \big|\, (\mathrm{p}_\ell S_{\iii})((0,1)) \cap (\mathrm{p}_\ell S_{\jjj})((0,1)) \neq \emptyset \big\} &\leq \sum_{x\in\{(\mathrm{p}_\ell S_{\jjj})(0),(\mathrm{p}_\ell S_{\jjj})(1/2),(\mathrm{p}_\ell S_{\jjj})(1)\}}\# \big\{\iii \in \Gamma_\ell' \, \big|\, x \in (\mathrm{p}_\ell S_{\iii})((0,1)) \big\}\\ &\leq 3(K \cdot g(r_{\min}^n))^m.
\end{align*}

For the second statement, we proceed similarly. 
By the definition \eqref{eq:projrow}, $\Sigma_{\jjj}'$ consists of $m$-tuples of symbols of functions with the same $\ell$-projection. Thus for any $\jjj \in \Sigma'$ and $\iii \in \Sigma_{\jjj}'$
\begin{align*}
S_{\iii}((0,1)^2) \cap S_{\jjj}((0,1)^2) \neq \emptyset \iff (\mathrm{p}_{3-\ell}S_{\iii})((0,1)) \cap (\mathrm{p}_{3-\ell}S_{\jjj})((0,1)) \neq \emptyset.
\end{align*}
Since $\Sigma'$ is \ref{B1}, for a fixed $\jjj \in \Sigma'$, if $x_k :=S_{\jjj}(y_k)$, $k = 1,2,3$, where $\mathrm{p}_\ell(y_k) =0$ and $\mathrm{p}_\ell(y_1) =0,\, \mathrm{p}_\ell(y_2) = \tfrac{1}{2},\, \mathrm{p}_\ell(y_3) = 1$, then for $\iii \in \Sigma_{\jjj}'$
\begin{align*}
S_{\iii}((0,1)^2) \cap S_{\jjj}((0,1)^2) \neq \emptyset \quad \Longrightarrow \quad \exists k \in \{1,2,3\}: x_k \in S_{\iii}((0,1)^2)
\end{align*}
Recall $\Delta_\delta$ from \eqref{eq:Delta}. Then 
\begin{align*}
\Sigma^{\{n\}} \subseteq \bigcup_{k = 0}^K \Delta_{r_{\max}^{n+k}} \quad \text{where} \quad K := \left\lceil n \Big( \frac{\log r_{\min}}{\log r_{\max}}-1\Big) \right\rceil.
\end{align*}
Since $\mathbb{G}$ satisfies the \ref{ANC} there is a subpolynomial function $f(\delta)$, such that for any $x \in \R^2$
\begin{align*}
\#\big\{\iii \in \Sigma'_{\jjj} \, \big|\, x \in S_{\iii}((0,1)^2) \big\} &\leq \#\big\{\iii \in \Sigma^{\{n\}m} \, \big|\, x \in S_{\iii}((0,1)^2) \big\} \\ &\leq \left(\sup_{x\in[0,1]^2}\#\big\{\iii \in \Sigma^{\{n\}} \, \big|\, x \in S_{\iii}((0,1)^2) \big\}\right)^m \leq (K \cdot f(r_{\min}^n))^m.
\end{align*}
Therefore, by a similar argument
\begin{align*}
\# \left\{\iii \in \Sigma_{\jjj}' \, \middle|\, S_{\iii}((0,1)^2) \cap S_{\jjj}((0,1)^2) \neq \emptyset \right\} \leq 9 (K \cdot f(r_{\min}^n))^m.
\end{align*}
\end{Proof}

\begin{lemma} \label{lemma:hdimroscsubsystem}
Let $n,m \in \mathbb{N}^+$ and let $\Sigma' \subseteq \Sigma^{\{n\}m}$ be such that $\{S_{\iii}\}_{\iii \in \Sigma'}$ defines a \ref{B1} \ref{A1} that satisfies the \ref{ANC} and the \ref{AWSC} for the $\ell$-projection, where $r_{\jjj,\ell} = r_\ell' \geq r_{3-\ell}' = r_{\jjj,3-\ell}$ for any $\jjj \in \Sigma'$. Then there exists $\Sigma^* \subseteq \Sigma'$ and a subexponential function $C(n)$ such that $\{S_{\iii}\}_{\iii \in \Sigma^*}$ defines a \ref{A1} satisfying the \ref{C1}, and its $\ell$-projection satisfies the OSC; moreover,
\begin{align*}
\frac{\log \bigg(\sum\limits_{\jjj\in\Gamma_{\ell}^{\{n\}m}}(\#\Sigma_{\jjj}^\ell )^{\frac{\log r_\ell'}{\log r_{3-\ell}'}}\bigg)}{-\log r_\ell'} 
\leq 
\frac{\log \bigg(\sum\limits_{\jjj\in\Gamma_{\ell}^{\{n\}m}} (\#\Sigma_{\jjj}^{\ell,*} )^{\frac{\log r_\ell'}{\log r_{3-\ell}'}}\bigg)}{-\log r_\ell'} 
+ \frac{1}{n}\log C(n),
\end{align*}
where $\Sigma_{\jjj}^{\ell,*} := \Sigma_{\jjj}^{\ell} \cap \Sigma^*$.
\end{lemma}

\begin{Proof}

For each $\jjj \in \Gamma_\ell'$, let $\Sigma_{\jjj}^{\ell,*} \subseteq \Sigma_{\jjj}^{\ell}$ be maximal such that $S_{\iii}((0,1)^2) \cap S_{\kkk}((0,1)^2) = \emptyset$ for all $\iii\neq \kkk, \in\Sigma_{\jjj}^{\ell,*}$.  Then by Lemma~\ref{lemma:WSCnew}\eqref{eq:WCSuse1}, there exists a subexponential function $C(n)$ such that
\[  \# \Sigma_{\jjj}^{\ell} \leq C(n)^{m} \cdot \# \Sigma_{\jjj}^{\ell,*}\]
Let $\Gamma_\ell^* \subseteq \Gamma_\ell'$ be such that 
\begin{itemize}
\item $S_{\iii}((0,1)^2) \cap S_{\kkk}((0,1)^2) = \emptyset$ for any $\iii, \kkk, \in \Gamma_\ell^*, \iii \neq \kkk$, 
\item for any $\iii \in \Gamma_\ell'$ either $\iii \in \Gamma_\ell^*$ or there is $\jjj \in \Gamma_\ell^*$ such that 
\[
(\mathrm{p}_\ell S_{\iii})((0,1)) \cap (\mathrm{p}_\ell S_{\jjj})((0,1)) \neq \emptyset \quad \text{ and } \quad \#\Sigma_{\jjj}^{\ell,*} \geq \#\Sigma_{\iii}^{\ell,*}.
\]
\end{itemize}
A simple greedy algorithm can define such a set. By Lemma~\ref{lemma:WSCnew}\eqref{eq:WSCuse2},
\[ \#\Gamma_\ell' \leq C(n)^{m} \cdot \#\Gamma_\ell^* .\]
Finally,
\begin{multline*}
\frac{\log \bigg(\sum\limits_{\jjj\in\Gamma_{\ell}^{\{n\}m}}(\#\Sigma_{\jjj}^\ell )^{\frac{\log r_\ell'}{\log r_{3-\ell}'}}\bigg)}{-\log r_\ell'} \leq  \frac{\log \bigg(C(n)^{m(1+\frac{\log r_\ell'}{\log r_{3-\ell}'})}\sum\limits_{\jjj\in\Gamma_{\ell}^*} (\#\Sigma_{\jjj}^{\ell,*} )^{\frac{\log r_\ell'}{\log r_{3-\ell}'}}\bigg)}{-\log r_\ell'} \\
\leq \frac{\log \bigg(\sum\limits_{\jjj\in\Gamma_{\ell}^*} (\#\Sigma_{\jjj}^{\ell,*} )^{\frac{\log r_\ell'}{\log r_{3-\ell}'}}\bigg)}{-\log r_\ell'} + \Big(1+\frac{\log r_{\min}}{\log r_{\max}}\Big)\frac{\log C(n)}{-n\log r_{\max}}.
\end{multline*}
Thus, $\Sigma^* := \cup_{\jjj \in \Gamma_\ell^*} \Sigma_{\jjj}^{\ell,*}$ gives the claimed subsystem.
\end{Proof}

We are now ready to prove Theorem \ref{T2}.

\subsection{Proof of Theorem \ref{T2}}

Let $\mathbb{G}$ be an \ref{B2} \ref{A1}, such that the principal projections satisfy the \ref{AWSC}. By Theorem \ref{thm:PS}, for any $n \in \N$
\begin{align*}
\dH(\Lambda) &\leq \mathbb{H}_{BA}(\mathbb{G}_n) = \max_{p \in P^{\{n\}}}D(p, \mathbb{G}_n)
\end{align*}
By Lemma \ref{lemma:symbolichomogene}, for any $\eps>0$ there is $M \in \N$ such that for any $m>M$ there is $\Sigma' \subseteq \Sigma^{\{n\}m}$ \ref{B1} \ref{A1} and $\ell \in \{1,2\}$ such that
\begin{align*}
\max_{p \in P^{\{n\}}}D(p, \mathbb{G}_n) \leq \frac{\log \bigg(\sum\limits_{\jjj\in\Gamma_{\ell}^{\{n\}m}}(\#\Sigma_{\jjj}^\ell )^{\frac{\log r_\ell'}{\log r_{3-\ell}'}}\bigg)}{-\log r_\ell'} + \eps.
\end{align*}
By Lemma \ref{lemma:hdimroscsubsystem}, there is $\Sigma^* \subseteq \Sigma'$ such that
\begin{align*}
\frac{\log \bigg(\sum\limits_{\jjj\in\Gamma_{\ell}^{\{n\}m}}(\#\Sigma_{\jjj}^\ell )^{\frac{\log r_\ell'}{\log r_{3-\ell}'}}\bigg)}{-\log r_\ell'} \leq \frac{\log \bigg(\sum\limits_{\jjj\in\Gamma_{\ell}^{*}} (\#\Sigma_{\jjj}^{\ell,*} )^{\frac{\log r_\ell'}{\log r_{3-\ell}'}}\bigg)}{-\log r_\ell'} + \frac{1}{n}\log C(n).
\end{align*}

Let $\Lambda_{n,m}^*$ be the attractor of the IFS $\{S_{\iii}\}_{\iii \in \Sigma^*}$. By Lemma \ref{lemma:hdimroscsubsystem}, $\Lambda_{n,m}^*$ satisfies the \ref{C1} and its $\ell$-projection satisfies the OSC. It inherits \ref{B2}, thus by Theorem \ref{thm:GL},
\begin{align*}
\frac{\log \bigg(\sum\limits_{\jjj \in \Gamma_{\ell}^{*}} (\#\Sigma_{\jjj}^{\ell,*})^{\frac{\log r_\ell'}{\log r_{3-\ell}'}}\bigg)}{-\log r_\ell'} = \dH(\Lambda_{n,m}^*).
\end{align*}
Therefore,
\begin{align*}
\dH(\Lambda) \leq \mathbb{H}_{BA}(\mathbb{G}_n) \leq \dH(\Lambda_{n,m}^*) + \frac{1}{n}\log C(n) + \eps \leq \dH(\Lambda) + \frac{1}{n}\log C(n) + \eps .
\end{align*}
Since $\varepsilon>0$ was arbitrary, one can finish the proof by taking $n\to\infty$.

\vspace{5pt}

\begin{remark}
In the case where $\mathbb{G}$ is a \ref{A1} satisfying \ref{G1}, the \ref{ANC}, and the \ref{AWSC} for the $y$-projection, the above proof works with the following modifications
\begin{itemize}
\item By the \ref{G1}, $\ell = 2$ in Lemma \ref{lemma:symbolichomogene}.
\item Thus, Lemma \ref{lemma:hdimroscsubsystem} is only needed and holds for $\ell = 2$.
\end{itemize}
\end{remark}

\vspace{5 pt}

\section{Box-counting dimension}

Now, we turn our attention to the box-counting dimension. Let us first state the corresponding theorem in a more general form, which clearly will imply Theorem~\ref{tb}.

\begin{theorem}[Box-counting dimension] \label{T1}
Let 
\[\mathbb{G} := \left\{S_i(x, y) := \left(r_{i,1}x+t_{i,1},r_{i,2}y+t_{i,2}\right)\right\}_{i\in\Sigma}
\]
be a \ref{A1} such that either the \ref{AWSC} holds for both projections, or it satisfies the \ref{G1}, the \ref{ANC} and the \ref{AWSC} for the $y$-projection. Then
\begin{align}\label{eq:firstbd}
\dB(\Lambda) = \max_{\ell \in \{1,2\}} \limsup_{\delta \to 0} \Bigg( \frac{\log \# \mathbb{S}_\delta^\ell}{-\log\delta}
+ \dB(\mathrm{p}_\ell\Lambda) \Bigg( 1+ \frac{\log \mathfrak{M}_{\dB(\mathrm{p}_\ell\Lambda)}\big( |r_{S,\ell}| \, \big|\, S \in \mathbb{S}_\delta^\ell\big)}{-\log\delta} \Bigg)\Bigg)
\end{align}
where $\mathfrak{M}_p(x_1, \dots , x_n ) := \Big(\frac{1}{n} \sum_{i=1}^nx_i^p \Big)^{\frac{1}{p}}$ is the power mean of the multiset $(x_1, \dots , x_n )$ with exponent $p$.
Moreover,
\begin{equation*}
 \dB(\Lambda) =\limsup_{\delta \to 0} \max\big\{ D_\delta^1, D_\delta^2\big\} = \limsup_{n \to \infty}\max\big\{d_n^1, d_n^2\big\},
\end{equation*}
where for $\ell \in \{1,2\}$, the quantities $D_\delta^\ell$ and $d_n^\ell$ are the unique real solutions of the equations
\begin{align*}
\sum_{S \in \mathbb{S}_\delta} \left( \frac{|r_{S,\ell}|}{|r_{S,3-\ell}|} \right)^{\dB(\mathrm{p}_\ell\Lambda)} |r_{S,3-\ell}|^{D_\delta^\ell} = 1 \quad \text{and} \quad
\sum_{S \in \mathbb{G}_n} \left( \frac{|r_{S,\ell}|}{|r_{S,3-\ell}|} \right)^{\dB(\mathrm{p}_\ell\Lambda)} |r_{S,3-\ell}|^{d_n^\ell} = 1.
\end{align*}
Furthermore, if $\mathrm{p}_\ell\Lambda$ satisfies the \ref{AWSC}, then instead of $D_\delta^\ell$, one can use $\mathfrak{D}_n^\ell$, the unique real solutions of the equations
\begin{align*}
\sum_{S \in \mathbb{G}_n} \left( \frac{|r_{S,\ell}|}{|r_{S,3-\ell}|} \right)^{s_n^\ell} |r_{S,3-\ell}|^{\mathfrak{D}_n^\ell} = 1,
\end{align*}
where $s_n^\ell$  is the similarity dimension of $\mathrm{p}_\ell\mathbb{G}_n$: 
\[\sum_{S\in\mathrm{p}_\ell\mathbb{G}_n}|r_{S,\ell}|^{s_n^\ell}=1.\]
\end{theorem}

For typical \ref{A1}s, the cardinality of the sets $\mathbb{S}_\delta^1$, $\mathbb{S}_\delta^2$ and $\mathbb{S}_\delta$ is much harder to compute compared to the $n$-th level functions. The remainder of the section is devoted to proving Theorem \ref{T1}.

\subsection{Proof of the initial formula in Theorem \ref{T1}}

For $\delta>0$, let $\Delta_\delta^1, \Delta_\delta^2$ be as in Subsection \ref{subsection:symbolic}. An equivalent characterisation of the box-counting dimension is obtained by defining $N_\delta$ via covers by axis-aligned squares of side length $\delta$. We adopt this formulation in the proof.

\begin{lemma} \label{NL3}
Let $\mathbb{G}$ be a \ref{A1}. Then for any $\varepsilon>0$,  there exists $c> 0$ such that for any $\iii \in \Delta_\delta^\ell$
\begin{equation*}\begin{aligned}
N_\delta(S_{\iii}(\Lambda)) \in \left(c^{-1} \left(\frac{|r_{\iii,\ell}|}{\delta}\right)^{\dB(\mathrm{p}_\ell\Lambda)-\varepsilon},\ c \left(\frac{|r_{\iii,\ell}|}{\delta}\right)^{\dB(\mathrm{p}_\ell\Lambda)+\varepsilon}\right)
\end{aligned}\end{equation*}
for every $\ell\in\{1,2\}$.
\end{lemma}
\begin{Proof} 

	Let $\ell \in \{1,2\}$. Observe that if $\iii \in \Delta_\delta^\ell$, then
	\begin{align*}
		N_\delta(S_{\iii}\big(\Lambda)\big) = N_{\delta \cdot |r_{\iii,\ell}|^{-1}}(\mathrm{p}_\ell\Lambda).
	\end{align*}
	See Figure \ref{Fig23} for an illustration.
	
	\begin{figure}[h]
		\centering
		\scalebox{.73}{
			\begin{tikzpicture}
				\draw (-2,0) rectangle (2,4);
				\draw[pattern=north west lines, pattern color=blue] (0.5,1) rectangle (0.7,3.4);
				\draw[pattern=north east lines, pattern color=blue] (0.5,1) rectangle (0.7,3.4);
				\draw[pattern=north west lines, pattern color=red] (0.35,1.1) rectangle (0.85,1.6);
				\draw[pattern=north east lines, pattern color=red] (0.35,1.1) rectangle (0.85,1.6);
				\draw[pattern=north west lines, pattern color=red] (0.35,2.1) rectangle (0.85,2.6);
				\draw[pattern=north east lines, pattern color=red] (0.35,2.1) rectangle (0.85,2.6);
				\draw[pattern=north west lines, pattern color=red] (0.35,2.5) rectangle (0.85,3);
				\draw[pattern=north east lines, pattern color=red] (0.35,2.5) rectangle (0.85,3);
				
				\draw[decorate,decoration={brace,amplitude=5pt,mirror,raise=2ex}] (-2,0) -- (2,0) node[midway,yshift=-2em]{1};
				\draw[decorate,decoration={brace,amplitude=5pt,raise=2ex}] (-2,0) -- (-2,4) node[midway,xshift=-2em]{1};
				\draw[decorate,decoration={brace,amplitude=5pt,raise=4ex}] (0.7,1) -- (0.7,3.4) node[midway,xshift=-3.5em]{$|r_{\iii,2}|$};
				\draw[decorate,decoration={brace,amplitude=5pt,raise=1ex, mirror}] (0.5,1) -- (0.7,1) node[midway,yshift=-1.5em]{$|r_{\iii,1}|$};
				\draw[decorate,decoration={brace,amplitude=5pt,raise=1ex, mirror}] (0.85,1.1) -- (0.85,1.6) node[midway,xshift=1.5em]{$\delta$};
				
				\draw[pattern={mylines[size= 6pt,line width=.6pt,angle=171]}, pattern color=blue] (8,0) rectangle (12,4);
				\draw[pattern={mylines[size= 20pt,line width=.6pt,angle=81]}, pattern color=blue] (8,0) rectangle (12,4);
				\draw[pattern={mylines[size= 6pt,line width=.6pt,angle=171]}, pattern color=red] (6,0.2) rectangle (14,1);
				\draw[pattern={mylines[size= 20pt,line width=.6pt,angle=81]}, pattern color=red] (6,0.2) rectangle (14,1);
				\draw[pattern={mylines[size= 6pt,line width=.6pt,angle=171]}, pattern color=red] (6,2.1) rectangle (14,2.9);
				\draw[pattern={mylines[size= 20pt,line width=.6pt,angle=81]}, pattern color=red] (6,2.1) rectangle (14,2.9);
				\draw[pattern={mylines[size= 6pt,line width=.6pt,angle=171]}, pattern color=red] (6,2.8) rectangle (14,3.6);
				\draw[pattern={mylines[size= 20pt,line width=.6pt,angle=81]}, pattern color=red] (6,2.8) rectangle (14,3.6);
				
				\draw[decorate,decoration={brace,amplitude=5pt,raise=2ex}] (8,4) -- (12,4) node[midway,yshift=2em]{1};
				\draw[decorate,decoration={brace,amplitude=5pt,mirror,raise=1ex}](6,0) -- (14,0) node[midway,yshift=-2em]{$\delta \cdot |r_{\iii,1}|^{-1}$};
				\draw[decorate,decoration={brace,amplitude=5pt,raise=1ex}] (6,0.2) -- (6,1) node[midway,xshift=-3em]{$\delta \cdot |r_{\iii,2}|^{-1}$};
				
				\draw[pattern={mylines[size= 6pt,line width=.6pt,angle=171]}, pattern color=blue] (8,-7) rectangle (12,-3);
				\draw[pattern={mylines[size= 20pt,line width=.6pt,angle=81]}, pattern color=blue] (8,-7) rectangle (12,-3);
				\draw[pattern={mylines[size= 6pt,line width=.6pt,angle=171]}, pattern color=red] (8,-6.8) rectangle (12,-6);
				\draw[pattern={mylines[size= 20pt,line width=.6pt,angle=81]}, pattern color=red] (8,-6.8) rectangle (12,-6);
				\draw[pattern={mylines[size= 6pt,line width=.6pt,angle=171]}, pattern color=red] (8,-4.9) rectangle (12,-4.1);
				\draw[pattern={mylines[size= 20pt,line width=.6pt,angle=81]}, pattern color=red] (8,-4.9) rectangle (12,-4.1);
				\draw[pattern={mylines[size= 6pt,line width=.6pt,angle=171]}, pattern color=red] (8,-4.2) rectangle (12,-3.4);
				\draw[pattern={mylines[size= 20pt,line width=.6pt,angle=81]}, pattern color=red] (8,-4.2) rectangle (12,-3.4);
				
				\draw[pattern={mylines[size= 6pt,line width=.6pt,angle=171]}, color=blue, thick] (0,-7) rectangle (0,-3);
				\draw[pattern={mylines[size= 6pt,line width=.6pt,angle=171]}, color=red, thick] (0.2,-6.8) rectangle (0.2,-6);
				\draw[pattern={mylines[size= 6pt,line width=.6pt,angle=171]}, color=red, thick] (0.16,-6.8) rectangle (0.24,-6.8);
				\draw[pattern={mylines[size= 6pt,line width=.6pt,angle=171]}, color=red, thick] (0.16,-6) rectangle (0.24,-6);
				
				\draw[pattern={mylines[size= 6pt,line width=.6pt,angle=171]}, color=red, thick] (0.2,-4.9) rectangle (0.2,-4.1);
				\draw[pattern={mylines[size= 6pt,line width=.6pt,angle=171]}, color=red, thick] (0.16,-4.9) rectangle (0.24,-4.9);
				\draw[pattern={mylines[size= 6pt,line width=.6pt,angle=171]}, color=red, thick] (0.16,-4.1) rectangle (0.24,-4.1);
				
				\draw[pattern={mylines[size= 6pt,line width=.6pt,angle=171]}, color=red, thick] (0.27,-4.2) rectangle (0.27,-3.4);
				\draw[pattern={mylines[size= 6pt,line width=.6pt,angle=171]}, color=red, thick] (0.23,-4.2) rectangle (0.31,-4.2);
				\draw[pattern={mylines[size= 6pt,line width=.6pt,angle=171]}, color=red, thick] (0.23,-3.4) rectangle (0.31,-3.4);
				
				\draw[decorate,decoration={brace,amplitude=5pt,raise=1ex, mirror}] (0.2,-6.8) -- (0.2,-6) node[midway,xshift=3.4em]{$\delta \cdot |r_{\iii,2}|^{-1}$};
				
				\draw [line width=1pt, double distance=3pt, arrows = {-Latex[length=0pt 3 0]}] (3,2) -- (5,2) ;
				\draw [line width=1pt, double distance=3pt, arrows = {-Latex[length=0pt 3 0]}] (10,-1.5) -- (10,-2.5);
				\draw [line width=1pt, double distance=3pt, arrows = {-Latex[length=0pt 3 0]}] (6,-5) -- (2,-5);
				
				\draw (4,2.4) node {$S_i^{-1}$};
				\draw (4,-4.6) node {$\mathrm{p}_2$};
				
				\draw[dashed,<-] (0.6,3.5)--(0.6,4.2) node[above,pos=1] {$S_{\iii}([0,1]^2)$};
				\draw[dashed,<-] (13,3.7)--(13,4.2) node[above,pos=1] {$S_{\iii}^{-1}(\delta \times \delta)$};
				\draw[dashed,<-] (12.1,-3.75)--(12.8,-3.75) node[right,pos=1] {$[0,1] \times (\mathrm{p}_2\circ S_{\iii})^{-1}(\delta)$};
		\end{tikzpicture}}\caption{\small The argument for a fixed cylinder and three $\delta$ by $\delta$ square transforming. Notice that for the cover of the cylinder by $\delta$ by $\delta$ squares, it is enough to consider the cover of the projection of the attractor to the $y$-axis by $\delta \cdot |r_{\iii,2}|^{-1}$ intervals.}
	\label{Fig23}
	\end{figure}
	
	\noindent For self-similar sets, and in particular for $\mathrm{p}_1\Lambda$ and $\mathrm{p}_2\Lambda$, the box-counting dimension exists, and therefore for any $\varepsilon>0$ there is $\Gamma > 0$ such that for any $\delta \leq \Gamma$
	\begin{align*}
		N_\delta(\mathrm{p}_\ell\Lambda) &\in \left(\delta^{-\dB(\mathrm{p}_\ell\Lambda)+\varepsilon},\delta^{-\dB(\mathrm{p}_\ell\Lambda)-\varepsilon}\right).
	\end{align*}
    Notice that $\delta \cdot |r_{\iii,\ell}|^{-1} \leq r_{\min}^{-1}$ for any $\iii \in \Delta_\delta^\ell$, and then let
		\begin{align*}
				\sup_{\iii \in \Delta_\delta^\ell} \left\{
				\frac{N_{\delta |r_{\iii,\ell}|^{-1}}(\mathrm{p}_\ell\Lambda)}
				{\left(\frac{\delta}{|r_{\iii,\ell}|}\right)^{-\dB(\mathrm{p}_\ell\Lambda)-\varepsilon}}
				\right\}
				&= \max \left\{ 1, \sup_{\substack{\iii \in \Delta_\delta^\ell \\ \delta |r_{\iii,\ell}|^{-1} > \Gamma}}
				\left\{ \frac{N_{\delta |r_{\iii,\ell}|^{-1}}(\mathrm{p}_\ell\Lambda)}
				{\left(\frac{\delta}{|r_{\iii,\ell}|}\right)^{-\dB(\mathrm{p}_\ell\Lambda)-\varepsilon}}
				\right\} \right\} \\
				&\leq \max \left\{ 1, N_{\Gamma}(\mathrm{p}_\ell\Lambda)
				\cdot r_{\min}^{-\dB(\mathrm{p}_\ell\Lambda)-    					\varepsilon}
				\right\} =: c_1
			\end{align*}
		and
		\begin{align*}
				\inf_{\iii \in \Delta_\delta^\ell} \left\{
				\frac{N_{\delta |r_{\iii,\ell}|^{-1}}(\mathrm{p}_\ell\Lambda)}
				{\left(\frac{\delta}{|r_{\iii,\ell}|}\right)^{-\dB(\mathrm{p}_\ell\Lambda)+\varepsilon}}
				\right\}
				&= \min \left\{ 1, \inf_{\substack{\iii \in \Delta_\delta^\ell \\ \delta |r_{\iii,\ell}|^{-1} > \Gamma}}
				\left\{ \frac{N_{\delta |r_{\iii,\ell}|^{-1}}(\mathrm{p}_\ell\Lambda)}
				{\left(\frac{\delta}{|r_{\iii,\ell}|}\right)^{-\dB(\mathrm{p}_\ell\Lambda)+\varepsilon}}
				\right\} \right\} \\
				&\geq \min \left\{ 1,
				N_{r_{\min}^{-1}}(\mathrm{p}_\ell\Lambda)
				\cdot \Gamma^{\dB(\mathrm{p}_\ell\Lambda)+\varepsilon}
				\right\} =: c_2.
			\end{align*}
Then the claim follows by the choice of the constant $c:= \max\{c_1, c_2^{-1}\}$. \end{Proof}
\begin{lemma}\label{lem:step1}
Let $\mathbb{G}$ be a \ref{A1}. Then 
\begin{equation*}
	\dUB(\Lambda)\leq\max_{\ell\in\{1,2\}}\limsup_{\delta\to0}\frac{\log\left(\sum_{S\in\mathbb{S}_{\delta}^{\ell}}\left(\frac{|r_{S,\ell}|}{\delta}\right)^{\dB(\mathrm{p}_\ell\Lambda)}\right)}{-\log\delta}.
\end{equation*}
\end{lemma}

\begin{Proof} 

It is easy to see that $\Lambda=\bigcup_{S \in \mathbb{S}_\delta}S(\Lambda)$ for every $\delta>0$, and clearly, $\mathbb{S}_\delta=\bigcup_{\ell\in\{1,2\}}\mathbb{S}_\delta^\ell$.
For every $\iii\in\Delta_\delta^\ell$,
\begin{equation*}
\delta r_{\min} < |r_{\iii,\ell}| \leq \delta^{\frac{\log(r_{\max})}{\log(r_{\min})}}.
\end{equation*}
Thus, by Lemma~\ref{NL3}, for every $\varepsilon>0$ there exists $c>0$ such that for every $\delta>0$
\begin{align*}
	N_\delta(\Lambda)&=N_{\delta}\left(\bigcup_{\ell\in\{1,2\}}\bigcup_{S\in\mathbb{S}_{\delta}^{\ell}}S(\Lambda)\right)\leq2\max_{\ell\in\{1,2\}}\sum_{S\in\mathbb{S}_{\delta}^{\ell}}N_{\delta}(S(\Lambda))\\
	&\leq 2c\max_{\ell\in\{1,2\}}\sum_{S\in\mathbb{S}_{\delta}^{\ell}}\left(\frac{|r_{S,\ell}|}{\delta}\right)^{\dB(\mathrm{p}_\ell\Lambda)+\varepsilon}\leq  2c\delta^{(\frac{\log(r_{\max})}{\log(r_{\min})}-1)\varepsilon}\max_{\ell\in\{1,2\}}\sum_{S\in\mathbb{S}_{\delta}^{\ell}}\left(\frac{|r_{S,\ell}|}{\delta}\right)^{\dB(\mathrm{p}_\ell\Lambda)}.
\end{align*}
Since $\varepsilon>0$ was arbitrary, Lemma \ref{lem:step1} follows.
\end{Proof}

\begin{theorem} \label{lem:step2}
Let $\mathbb{G}$ be an \ref{A1} satisfying the \ref{AWSC} for both principal projections. Then
\begin{equation*}
\max_{\ell\in\{1,2\}}\limsup_{\delta\to0}\frac{\log\left(\sum_{S\in\mathbb{S}_{\delta}^{\ell}}\left(\frac{|r_{S,\ell}|}{\delta}\right)^{\dB(\mathrm{p}_\ell\Lambda)}\right)}{-\log\delta}\leq\dLB(\Lambda).
\end{equation*}
\end{theorem}

The proof will have a similar flavour to the proof for the Hausdorff dimension. We will construct a sufficiently well-separated subsystem: By allowing an $\eps$ error, we will decompose $\mathbb{S}_\delta^\ell$ to a finite number of subsystems with close to \ref{B1} contractions. Then by the \ref{AWSC}, we will be able to find a subsystem with close dimension while satisfying the \ref{C1} and the OSC for the $\ell$-projection. By enlarging this IFS with functions from $\mathbb{S}^{3-\ell}$, we will grow the dimension of the $\ell$-projection, while not compromising the \ref{C1} and the OSC for the $\ell$-projection.

\vspace{5pt}

\begin{Proof}

Let $\ell\in\{1,2\}$, $\eps > 0$ and $\delta > 0$. Decompose $\mathbb{S}_\delta^\ell$ with respect to its larger contraction ratio: for every $S\in\mathbb{S}_\delta^\ell$, $\delta r_{\min} < |r_{S,\ell}|\leq \delta^{\frac{\log r_{\max}}{\log r_{\min}}}$. Thus,
\begin{equation*}
\mathbb{S}_\delta^\ell \subseteq \bigcup_{n=\left\lfloor\frac{\log r_{\max}}{\varepsilon\log r_{\min}}\right\rfloor}^{\lceil1/\varepsilon\rceil}\mathbb{S}_{\delta, \varepsilon}^{\ell, n}, \quad \text{where} \quad \mathbb{S}_{\delta, \varepsilon}^{\ell, n} := \Big\{S \in \mathbb{S}_\delta^\ell \, \Big|\, \delta^{(n+1)\varepsilon}\leq|r_{S,\ell}|\leq\delta^{n\varepsilon} \Big\}.
\end{equation*}
We now construct a well-separated subsystem $\mathbb{S}_{\delta, \varepsilon}^{\ell, n, *} \subseteq \mathbb{S}_{\delta, \varepsilon}^{\ell, n}$. Let us first decompose $\mathbb{S}_{\delta, \varepsilon}^{\ell, n}$ according to the maps with common principal projection. That is, for a $T \in \mathrm{p}_\ell \mathbb{S}_{\delta, \varepsilon}^{\ell, n}$, let
\begin{align*}
	 \quad R(T) := \left\{ S \in \mathbb{S}_{\delta, \varepsilon}^{\ell, n} \ \middle| \ \mathrm{p}_\ell S = T\right\}.
\end{align*}
Then clearly,
\begin{align*}
    \mathbb{S}_{\delta, \varepsilon}^{\ell, n} = \bigcup_{T \in \mathrm{p}_\ell \mathbb{S}_{\delta, \varepsilon}^{\ell, n}} R(T).
\end{align*}
For every $T \in \mathrm{p}_\ell \mathbb{S}_{\delta, \varepsilon}^{\ell, n}$, let $R^*(T) \subseteq R(T)$ be a maximal subset with respect to cardinality such that $S((0,1)^2)\cap \hat{S}((0,1)^2)=\emptyset$ for every $S\neq\hat{S}\in R^*(T)$. Then let $\mathrm{p}_\ell\mathbb{S}_{\delta, \varepsilon}^{\ell, n,*} \subseteq \mathrm{p}_\ell \mathbb{S}_{\delta, \varepsilon}^{\ell, n}$ be such that
\begin{itemize}
\item $T((0,1)) \cap \hat{T}((0,1)) = \emptyset$ for every $T \neq \hat{T}$, $T, \hat{T} \in \mathrm{p}_\ell\mathbb{S}_{\delta, \varepsilon}^{\ell, n,*}$,
\item for every $T \in \mathrm{p}_\ell\mathbb{S}_{\delta, \varepsilon}^{\ell, n}$ either $T \in \mathrm{p}_\ell\mathbb{S}_{\delta, \varepsilon}^{\ell, n,*}$ or there is another $\hat{T} \in \mathrm{p}_\ell\mathbb{S}_{\delta, \varepsilon}^{\ell, n,*}$ such that $\#R^*(T) \leq \#R^*(\hat{T})$ and $T((0,1))\cap \hat{T}((0,1)) \neq \emptyset$.
\end{itemize}
Such a set can be constructed clearly by a simple greedy algorithm. Finally, let
\begin{equation*}
\mathbb{S}_{\delta, \varepsilon}^{\ell, n, *} := \bigcup_{T \in \mathrm{p}_\ell\mathbb{S}_{\delta, \varepsilon}^{\ell, n, *}} R^*(T).
\end{equation*}

\begin{lemma} \label{futol1} There exists a subpolynomial function $g(\delta)$ such that 
	\begin{align*}
		\# \mathbb{S}_{\delta, \varepsilon}^{\ell, n} \leq g(\delta) \cdot \delta^{-\eps} \cdot \# \mathbb{S}_{\delta, \varepsilon}^{\ell, n, *}.
	\end{align*}
\end{lemma}

\begin{Proof}

We estimate each projection separately, and then the product of the two bounds gives the claimed bound.
By definition, 
\[
R(T) \subseteq \mathrm{p}_{3-\ell}\mathbb{S}_{\delta, \eps}^{\ell, n} \subseteq \mathrm{p}_{3-\ell}\mathbb{S}_{\delta}^{\ell} \subseteq \mathcal{M}_\delta^{3-\ell}\text{ for $T \in \mathrm{p}_{3-\ell}\mathbb{S}_{\delta, \eps}^{\ell, n}$}.
\]
Since $\mathrm{p}_{3-\ell}\mathbb{G}$ satisfies the \ref{AWSC}, there is a subpolynomial function $g_1(\delta)$ such that 
\begin{align*}
\#R(T) \leq g_1(\delta) \cdot \#R^*(T).
\end{align*}

On the other hand, by definition
\begin{align*}
\mathrm{p}_{\ell}\mathbb{S}_{\delta, \eps}^{\ell, n} \subseteq \mathrm{p}_{\ell}\big\{S_{\iii}\, \big|\, \delta^{(n+1)\varepsilon}\leq|r_{S,\ell}|\leq\delta^{n\varepsilon} \big\} \subseteq \bigcup_{m = 0}^{A}\mathcal{M}_{\delta^{n\eps}r_{\max}^m}^\ell,
\end{align*}
where $A = \big\lceil\eps\frac{\log \delta}{\log r_{\max}}\big\rceil$. For any $S \in \mathbb{S}_{\delta, \eps}^{\ell, n}$
\begin{itemize}
\item $\delta^{(n+1)\varepsilon} \leq|r_{S,\ell}|$, so a $\delta^{(n+1)\varepsilon}$-grid intersects each $(\mathrm{p}_\ell S)((0,1))$,
\item $|r_{S,\ell}|\leq\delta^{n\varepsilon}$, hence $(\mathrm{p}_\ell S)((0,1))$ contains at most $\lceil\frac{\delta^{n\eps}}{\delta^{(n+1)\eps}}\rceil+1 = \lceil\delta^{-\eps}\rceil+1$ points from the $\delta^{(n+1)\varepsilon}$-grid.
\end{itemize}
Additionally, since $\mathrm{p}_{\ell}\mathbb{G}$ satisfies the \ref{AWSC},
\begin{itemize}
\item there is a subpolynomial function $g_2(\delta)$ such that for each $x\in \R$, and in particular, for the points in the $\delta^{(n+1)\eps}$-grid, at most $g_2(\delta)$ cylinders contain $x$ from any fixed $\mathcal{M}_{\delta^{n\eps}r_{\max}^m}^\ell$.
\end{itemize}
Combining the above three estimates,
\begin{align*}
\#\big\{\mathrm{p}_\ell S_{\jjj} \in \mathrm{p}_{\ell}\mathbb{S}_{\delta, \eps}^{\ell, n}\, \big|\, (\mathrm{p}_\ell S_{\iii})((0,1)) \cap (\mathrm{p}_\ell S_{\jjj})((0,1)) \neq \emptyset \big\} \leq A \cdot g_2(\delta) \cdot (\lceil\delta^{-\eps}\rceil+1) \leq g_3(\delta) \cdot \delta^{-\eps}
\end{align*}
where $g_3(\delta)$ is a subpolynomial function. One combines the bounds for the two projections to have the claimed bound in the lemma.
\end{Proof}

Given Lemma \ref{futol1}, the approximation of the box dimension only follows if the constructed subsystem approximates the dimension of the projection, too; however, this is not necessarily the case. Hence, we have to enlarge $\mathbb{S}_{\delta, \varepsilon}^{\ell, n, *}$. Let $\mathbb{F} \subseteq \{S \, |\, \mathrm{p}_\ell S \in \mathcal{M}_{\delta^{n\varepsilon}}^{\ell}\}$ be a maximal subset with respect to containment such that the images of the unit interval with respect to the maps $\mathrm{p}_\ell\mathbb{F}$ are pairwise disjoint and for every $T\in\mathbb{F}$ and $S\in\mathbb{S}_{\delta, \varepsilon}^{\ell, n, *}$, $(\mathrm{p}_\ell T)((0,1))\cap(\mathrm{p}_\ell S)((0,1))=\emptyset$. Such a set can be defined again by a greedy algorithm. Denote by $\Lambda_{\delta,\varepsilon}^{\ell,n}$ the attractor of the IFS $\mathbb{F} \cup \mathbb{S}_{\delta,\epsilon}^{\ell,n,*}$. By construction, $\Lambda_{\delta,\varepsilon}^{\ell,n}$ satisfies the \ref{C1} and $\mathrm{p}_\ell \Lambda_{\delta,\varepsilon}^{\ell,n}$ satisfies the OSC. These two properties lead to the following Lemmas \ref{lemma:dimproj} and \ref{lemma:dim2d}.

\begin{lemma}\label{lemma:dimproj} There is a constant $c> 0$ and a subpolynomial function $g(\delta)$ such that 
\begin{align*}
\dB(\mathrm{p}_\ell\Lambda) &\leq \dB(\mathrm{p}_\ell \Lambda_{\delta,\varepsilon}^{\ell,n}) + \frac{\log g(\delta)}{-\log\delta}+ c\varepsilon.
\end{align*}
\end{lemma}
\begin{Proof}

Since $\mathrm{p}_\ell \Lambda_{\delta,\varepsilon}^{\ell,n}$ satisfies the OSC, we get that
$$
\sum_{T\in\mathrm{p}_\ell \mathbb{F} \cup \mathrm{p}_\ell \mathbb{S}_{\delta,\epsilon}^{\ell,n,*}}|r_{T,\ell}|^{\dB(\mathrm{p}_\ell \Lambda_{\delta,\varepsilon}^{\ell,n})}=1.
$$
By definition, the cover
$$
\{T((0,1))\, \big|\, T\in\mathrm{p}_\ell \mathbb{F} \cup \mathrm{p}_\ell \mathbb{S}_{\delta,\epsilon}^{\ell,n,*}\}
$$
forms a $\delta^{(n+1)\varepsilon}$-packing and a $3\delta^{n\varepsilon}$-cover of $\mathrm{p}_\ell\Lambda$. Hence, for a $g(\delta)$ subpolynomial,
\begin{align*}
1&\geq \delta^{(n+1)\varepsilon\dB(\mathrm{p}_\ell \Lambda_{\delta,\varepsilon}^{\ell,n})}\#(\mathrm{p}_\ell \mathbb{F} \cup \mathrm{p}_\ell \mathbb{S}_{\delta,\epsilon}^{\ell,n,*})\\
&\geq g(\delta)^{-1} \delta^{(n+1)\varepsilon\dB(\mathrm{p}_\ell \Lambda_{\delta,\varepsilon}^{\ell,n})}\delta^{-n\varepsilon(\dB(\mathrm{p}_\ell\Lambda)-\varepsilon)}\\
&\geq g(\delta)^{-1} \delta^{n\varepsilon(\dB(\mathrm{p}_\ell \Lambda_{\delta,\varepsilon}^{\ell,n})-\dB(\mathrm{p}_\ell\Lambda)+\varepsilon)+\varepsilon}.
\end{align*}
Finally, since $n\geq\frac{\log r_{\max}}{2\varepsilon\log r_{\min}}$, 
\begin{align*}
\dB(\mathrm{p}_\ell\Lambda)&\leq \dB(\mathrm{p}_\ell \Lambda_{\delta,\varepsilon}^{\ell,n}) + \frac{1}{n\eps}\frac{\log g(\delta)}{-\log \delta} - \eps +\frac{1}{n}\\
&\leq \dB(\mathrm{p}_\ell \Lambda_{\delta,\varepsilon}^{\ell,n}) + \frac{2\log r_{\min}}{\log r_{\max}} \frac{\log g(\delta)}{-\log\delta} + \left(\frac{2\log r_{\min}}{\log r_{\max}}-1\right)\varepsilon.
\end{align*}
\end{Proof}

\begin{lemma} \label{lemma:dim2d}
Let $d \in \R$ be the unique solution of the equation
\begin{equation*}
\sum_{S\in \mathbb{F} \cup \mathbb{S}_{\delta,\epsilon}^{\ell,n,*}}\left(\frac{|r_{S,\ell}|}{|r_{S,3-\ell}|}\right)^{\dB(\mathrm{p}_\ell \Lambda_{\delta,\varepsilon}^{\ell,n})}|r_{S,3-\ell}|^d=1.
\end{equation*}
Then $\dB(\Lambda_{\delta,\varepsilon}^{\ell,n})\geq d$.
\end{lemma}

\begin{Proof}
Since the IFS $\mathbb{F} \cup \mathbb{S}_{\delta,\epsilon}^{\ell,n,*}$ satisfies the \ref{C1} by construction, it follows from Theorem~\ref{thm:FW}.
\end{Proof}

Now we are ready to finish the proof of Theorem \ref{lem:step2}.
Let $\ell\in\{1,2\}$ and $\eps > 0$ be arbitrary. Let $N=\left\{\left\lfloor\frac{\log r_{\max}}{\varepsilon\log r_{\min}}\right\rfloor,\ldots,\lceil1/\varepsilon\rceil\right\}$. By the definition of $\mathbb{S}_{\delta,\epsilon}^{\ell,n}$,
\begin{align*}
\sum_{S\in\mathbb{S}_{\delta}^{\ell}}\left(\frac{|r_{S,\ell}|}{\delta}\right)^{\dB(\mathrm{p}_\ell\Lambda)}&\leq c_1 \max_{n\in N}\left\{\#\mathbb{S}_{\delta,\epsilon}^{\ell,n}\cdot\delta^{\min\{(n\varepsilon-1),0\}\dB(\mathrm{p}_\ell\Lambda)}\right\}.
\intertext{By Lemma \ref{futol1}, there is a subpolynomial function $g(\delta)$ such that}
&\leq c_2 g(\delta) \delta^{- \varepsilon}\max_{n\in N}\left\{\#\mathbb{S}_{\delta,\epsilon}^{\ell,n,*}\cdot\delta^{(n\varepsilon-1)\dB(\mathrm{p}_\ell\Lambda)}\right\}.
\intertext{By Lemma \ref{lemma:dimproj}, $\delta^{(n\varepsilon-1)\dB(\mathrm{p}_\ell\Lambda)} \leq \max\{g_1(\delta)^{(1-n\eps)}, 1\} \cdot \delta^{\min\{\eps(n\eps-1) C,0\}} \cdot \delta^{(n\varepsilon-1)\dB(\mathrm{p}_\ell\Lambda_{\delta,\varepsilon}^{\ell,n})}$, so}
&\leq c_3 g_2(\delta) \delta^{-c_4 \cdot \varepsilon}\max_{n\in N}\left\{\#\mathbb{S}_{\delta,\epsilon}^{\ell,n,*}\cdot\delta^{(n\varepsilon-1)\dB(\mathrm{p}_\ell \Lambda_{\delta,\varepsilon}^{\ell,n})}\right\} .
\intertext{By $\frac{|r_{S,\ell}|}{|r_{S,3-\ell}|} \geq \min\{1, r_{\min}\cdot\delta^{(n+1)\eps-1}\}$,}
&\leq c_3 g_2(\delta) \delta^{-c_5 \cdot \varepsilon}\max_{n\in N}\left\{\sum_{S\in\mathbb{S}_{\delta,\epsilon}^{\ell,n,*}}\left(\frac{|r_{S,\ell}|}{|r_{S,3-\ell}|}\right)^{\dB(\mathrm{p}_\ell \Lambda_{\delta,\varepsilon}^{\ell,n})}\right\}.
\intertext{By $|r_{S,3-\ell}| \geq r_{\min}\cdot \delta$, thus}
&\leq c_6 g_2(\delta) \delta^{-c_5\cdot \varepsilon -d}\max_{n\in N}\left\{\sum_{S\in\mathbb{S}_{\delta,\epsilon}^{\ell,n,*}}\left(\frac{|r_{S,\ell}|}{|r_{S,3-\ell}|}\right)^{\dB(\mathrm{p}_\ell \Lambda_{\delta,\varepsilon}^{\ell,n})}|r_{S,3-\ell}|^d \right\}.
\intertext{By Lemma \ref{lemma:dim2d}, $d\leq\dB(\Lambda_{\delta,\varepsilon}^{\ell,n})\leq\dLB(\Lambda)$, and so,}
&= c_6 g_2(\delta) \delta^{-c_5\cdot \varepsilon-d} \leq c_6 g_2(\delta) \delta^{-c_5\cdot \varepsilon-\dB(\Lambda_{\delta,\varepsilon}^{\ell,n})} \leq c_6 g_2(\delta) \delta^{-(c_5\cdot \varepsilon+\dLB(\Lambda))}.
\end{align*}
Since $c_1,\dots,c_6$ are constants independent of $\delta$ and $g_2(\delta)$ is subpolynomial in terms of $\delta$,
\[
\limsup_{\delta\to0}\frac{\log\left(\sum_{S\in\mathbb{S}_{\delta}^{\ell}}\left(\frac{|r_{S,\ell}|}{\delta}\right)^{\dB(\mathrm{p}_\ell\Lambda)}\right)}{-\log\delta}\leq\dLB(\Lambda)+c_5\cdot\eps.
\]
Since $n\in N$, $c_5 = 2 + \max\{0,(n\eps-1)\} \leq 3$, and since $\varepsilon>0$ was arbitrarily small, the claim follows.
\end{Proof}

\begin{theorem} \label{lem:step2nbh}
Let $\mathbb{G}$ be a \ref{A1} satisfying \ref{G1}, the \ref{ANC} and the \ref{AWSC} for the $y$-projection. Then
\begin{equation*}
\limsup_{\delta\to0}\frac{\log\left(\sum_{S\in\mathbb{S}_{\delta}^{2}}\left(\frac{|r_{S,2}|}{\delta}\right)^{\dB(\mathrm{p}_2\Lambda)}\right)}{-\log\delta}\leq\dLB(\Lambda).
\end{equation*}
\end{theorem}

\begin{Proof}
One mirrors the arguments of the proof of Theorem~\ref{lem:step2} with $\ell = 2$. The only part that changes is a small part in proof of Lemma \ref{futol1}. Namely,
\begin{align*}
\#R(T) \leq g_1(\delta) \cdot \#R^*(T)
\end{align*} 
for a $g_1(\delta)$ subpolynomial function, does not follow from the \ref{AWSC} for $\mathrm{p}_{3-\ell}\mathbb{G}$, but from the \ref{ANC}.
\end{Proof}

\begin{corollary}\label{cor:dim1}
	Let $\mathbb{G}$ be an \ref{A1} such that either the \ref{AWSC} holds for both projections, or it satisfies the \ref{G1}, the \ref{ANC} and the \ref{AWSC} for the $y$-projection. Then the box-counting dimension of the attractor exists. Moreover,
	\begin{equation*}
		\dB(\Lambda)=\max_{\ell\in\{1,2\}}\limsup_{\delta\to0}\frac{\log\left(\sum_{S\in\mathbb{S}_{\delta}^{\ell}}\left(\frac{|r_{S,\ell}|}{\delta}\right)^{\dB(\mathrm{p}_\ell\Lambda)}\right)}{-\log\delta}.
	\end{equation*}
\end{corollary}

\begin{Proof}
	It is an immediate corollary of Lemma \ref{lem:step1}, Theorem \ref{lem:step2nbh} and Theorem~\ref{lem:step2}.
\end{Proof}

A straightforward calculation shows that
\begin{align*}
\frac{\log \left(\sum_{S \in \mathbb{S}_\delta^\ell} \left(\frac{|r_{S,\ell}|}{\delta}\right)^{\dB(\mathrm{p}_\ell\Lambda)}\right)}{-\log\delta}
&= \frac{\log \# \mathbb{S}_\delta^\ell}{-\log\delta}
+ \dB(\mathrm{p}_\ell\Lambda) \left( 1 + \frac{\log\mathfrak{M}_{p}\left( |r_{S,\ell}| \,\middle|\, S \in \mathbb{S}_\delta^\ell\right)}{-\log\delta} \right),
\end{align*}
where $\mathfrak{M}_p$ denotes the power mean $\mathfrak{M}_p(x_1, \dots , x_n) = \Big(\frac{1}{n} \sum_{i=1}^nx_i^p \Big)^{\frac{1}{p}}$ with exponent $p = \dB(\mathrm{p}_\ell\Lambda) \in [0,1]$ and with $n = \#\mathbb{S}_\delta^\ell$. This completes the proof of \eqref{eq:firstbd}.

\subsection{Proof of the secondary formulas in Theorem \ref{T1}} \label{subsec:secform}
For $\ell \in \{1,2\}$, let us define $D_\delta^\ell$ by the equations:
\begin{align*}
\sum_{S \in \mathbb{S}_\delta}\left(\frac{|r_{S,\ell}|}{|r_{S,3-\ell}|}\right)^{\dB(\mathrm{p}_\ell\Lambda)}|r_{S,3-\ell}|^{D_\delta^\ell} = 1
\end{align*}
Denote $D_*^\ell := \limsup_{\delta \to 0} D_\delta^\ell$. Notice that by Corollary~\ref{cor:dim1}
\begin{align}
\dB(\Lambda) &= \max_{\ell \in \{1,2\}} \limsup_{\delta \to 0} \frac{\log \left(\sum_{S \in \mathbb{S}_\delta^\ell} \left(\frac{|r_{S,\ell}|}{|r_{S,3-\ell}|}\right)^{\dB(\mathrm{p}_\ell\Lambda)}\right)}{-\log\delta} \notag \\
&= \limsup_{\delta \to 0} \frac{\log \left(\sum_{\ell \in \{1,2\}} \delta^{-D_\delta^\ell}\sum_{S \in \mathbb{S}_\delta^\ell} \left(\frac{|r_{S,\ell}|}{|r_{S,3-\ell}|}\right)^{\dB(\mathrm{p}_\ell\Lambda)} |r_{S,3-\ell}|^{D_\delta^\ell} \right)}{-\log\delta} \notag \\
&\leq \limsup_{\delta \to 0} \frac{\log \left(\sum_{\ell \in \{1,2\}} \delta^{-D_\delta^\ell}\sum_{S \in \mathbb{S}_\delta} \left(\frac{|r_{S,\ell}|}{|r_{S,3-\ell}|}\right)^{\dB(\mathrm{p}_\ell\Lambda)} |r_{S,3-\ell}|^{D_\delta^\ell} \right)}{-\log\delta} \notag \\
&= \limsup_{\delta \to 0} \frac{\log \left(\sum_{\ell \in \{1,2\}} \delta^{-D_\delta^\ell}\right)}{-\log\delta} = \max_{\ell \in \{1,2\}} \left\{\limsup_{\delta \to 0} D_\delta^\ell\right\}.\label{eq:2ub}
\end{align}
For the other direction, let $\ell_\delta^*$ be such that $D_\delta^{\ell_\delta^*}= \max_{\ell \in \{1,2\}} \left\{D_\delta^\ell\right\}$. The following two lemmas are similar to the proof of Theorem~11.4.2 in the book of Bárány, Simon, and Solomyak \cite{MR4661364}.

\begin{lemma}\label{lem:tech}
	$\dB(\mathrm{p}_1\Lambda)+\dB(\mathrm{p}_2\Lambda)\geq \limsup_{\delta \to 0}D_\delta^{\ell_\delta^*}$.
\end{lemma}
\begin{Proof} 

Firstly, by the definition of $D_\delta^{\ell_\delta^*}$,
	\begin{align*}
		1 = \sum_{S \in \mathbb{S}_\delta} |r_{S,{\ell_\delta^*}}|^{\dB(\mathrm{p}_{\ell_\delta^*}(\Lambda))} |r_{S,3-{\ell_\delta^*}}|^{-\dB(\mathrm{p}_{\ell_\delta^*}(\Lambda))+D_\delta^{\ell_\delta^*}}.
	\end{align*}
Since $|r_{S,\ell}|\in\left[r_{\min}\delta,\delta^{\frac{\log r_{\max}}{\log r_{\min}}}\right]$ for $\ell \in \{1,2\}$ and $S\in\mathbb{S}_{\delta}$, it is enough to see that
	\begin{align} \label{eq:bpn1}
		0 \geq \limsup_{\delta\to0}\frac{\log\left( \sum_{S \in \mathbb{S}_\delta} |r_{S,1}|^{\dB(\mathrm{p}_1\Lambda)} |r_{S,2}|^{\dB(\mathrm{p}_2\Lambda)} \right)}{-\log\delta}.
	\end{align}
	Let us argue by contradiction. Suppose the opposite. Now
	\begin{align*}
		\dB(\Lambda) &= \limsup_{\delta \to 0} \frac{\log \left(\sum_{\ell \in \{1,2\}} \sum_{S \in \mathbb{S}_\delta^\ell} \left(\frac{|r_{S,\ell}|}{|r_{S,3-\ell}|}\right)^{\dB(\mathrm{p}_\ell\Lambda)}\right)}{-\log\delta} \\
		&= \limsup_{\delta \to 0} \frac{\log \left(\delta^{-\dB(\mathrm{p}_1\Lambda) - \dB(\mathrm{p}_2\Lambda)}\sum_{S \in \mathbb{S}_\delta} \left(|r_{S,1}|^{\dB(\mathrm{p}_1\Lambda)} |r_{S,2}|^{\dB(\mathrm{p}_2\Lambda)}\right)\right)}{-\log \delta}\\
		&= \dB(\mathrm{p}_1\Lambda) + \dB(\mathrm{p}_2\Lambda) + \limsup_{\delta \to 0} \frac{\log \left(\sum_{S \in \mathbb{S}_\delta} |r_{S,1}|^{\dB(\mathrm{p}_1\Lambda)} |r_{S,2}|^{\dB(\mathrm{p}_2\Lambda)}\right)}{-\log \delta}.
\end{align*}
This implies
\[ \dB(\Lambda) > \dB(\mathrm{p}_1\Lambda) + \dB(\mathrm{p}_2\Lambda),\]
which is a contradiction, since $\Lambda \subseteq \mathrm{p}_1\Lambda \times \mathrm{p}_2\Lambda$, and $
\dB(A)+\dB(B) \geq \dUB(A\times B)$ for any two sets, $A, B \subseteq \R^d$. See Falconer’s book \cite{Falconer1990}.
\end{Proof}
	
\begin{lemma} \label{Lemma:bpv1} There is a constant $c > 0$ such that
\begin{align*}
\limsup_{\delta \to 0} \max_{S \in \mathbb{S}_\delta^{3-\ell_\delta^*}} \frac{\left(\frac{|r_{S,3-\ell_\delta^*}|}{|r_{S,\ell_\delta^*}|}\right)^{\dB(\mathrm{p}_{3-\ell_\delta^*}(\Lambda))}|r_{S,\ell_\delta^*}|^{D_\delta^{3-\ell_\delta^*}}}{\left(\frac{|r_{S,\ell_\delta^*}|}{|r_{S,3-\ell_\delta^*}|}\right)^{\dB(\mathrm{p}_{\ell_\delta^*}(\Lambda))}|r_{S,3-\ell_\delta^*}|^{D_\delta^{\ell_\delta^*}}\delta^{D_\delta^{3-\ell_\delta^*}-D_\delta^{\ell_\delta^*}}} \geq c.
\end{align*}
\end{lemma}
\begin{Proof}
Let $S \in \mathbb{S}_\delta^{3-\ell_\delta^*}$. Observe that
\begin{align*}
 \frac{\left(\frac{|r_{S,3-\ell_\delta^*}|}{|r_{S,\ell_\delta^*}|}\right)^{\dB(\mathrm{p}_{3-\ell_\delta^*}(\Lambda))}|r_{S,\ell_\delta^*}|^{D_\delta^{3-\ell_\delta^*}}}{\left(\frac{|r_{S,\ell_\delta^*}|}{|r_{S,3-\ell_\delta^*}|}\right)^{\dB(\mathrm{p}_{\ell_\delta^*}(\Lambda))}|r_{S,3-\ell_\delta^*}|^{D_\delta^{\ell_\delta^*}}} = \left(\frac{|r_{S,3-\ell_\delta^*}|}{|r_{S,\ell_\delta^*}|}\right)^{\dB(\mathrm{p}_{3-\ell_\delta^*}(\Lambda))+\dB(\mathrm{p}_{\ell_\delta^*}(\Lambda))-D_\delta^{\ell_\delta^*}}
\cdot |r_{S,\ell_\delta^*}|^{D_\delta^{3-\ell_\delta^*}-D_\delta^{\ell_\delta^*}}.
\end{align*}
Since $\frac{|r_{S,3-\ell_\delta^*}|}{|r_{S,\ell_\delta^*}|} \geq 1$ and $\delta \geq |r_{S,\ell_\delta^*}| > \delta \cdot r_{\min}$, the claim follows by Lemma~\ref{lem:tech}.
\end{Proof}

\noindent Using Lemma \ref{Lemma:bpv1}, we get
\begin{align*}
\dB(\Lambda) &= \limsup_{\delta \to 0} \frac{\log \left(\sum_{\ell \in \{1,2\}} \sum_{S \in \mathbb{S}_\delta^\ell} \left(\frac{|r_{S,\ell}|}{|r_{S,3-\ell}|}\right)^{\dB(\mathrm{p}_\ell\Lambda)}\right)}{-\log\delta} \\
&\geq \limsup_{\delta \to 0} \frac{1}{-\log\delta} \Bigg(\log \Bigg(\delta^{-D_\delta^{\ell_\delta^*}}\sum_{S \in \mathbb{S}_\delta^{\ell_\delta^*}} \left(\frac{|r_{S,\ell_\delta^*}|}{|r_{S,3-\ell_\delta^*}|}\right)^{\dB(\mathrm{p}_{\ell_\delta^*}(\Lambda))}|r_{S,3-\ell_\delta^*}|^{-D_\delta^{\ell_\delta^*}} \\
& \quad \quad \quad \quad \quad  + \delta^{-D_\delta^{3-\ell_\delta^*}}\sum_{S \in \mathbb{S}_\delta^{3-\ell_\delta^*}} c \delta^{D_\delta^{3-\ell_\delta^*}-D_\delta^{\ell_\delta^*}} \left(\frac{|r_{S,\ell_\delta^*}|}{|r_{S,3-\ell_\delta^*}|}\right)^{\dB(\mathrm{p}_{\ell_\delta^*}(\Lambda))}|r_{S,3-\ell_\delta^*}|^{D_\delta^{\ell_\delta^*}} \Bigg) \Bigg) \\
&=\limsup_{\delta \to 0} \frac{\log \left(\delta^{-D_\delta^{\ell_\delta^*}}\sum_{S \in \mathbb{S}_\delta} \left(\frac{|r_{S,\ell_\delta^*}|}{|r_{S,3-\ell_\delta^*}|}\right)^{\dB(\mathrm{p}_{\ell_\delta^*}(\Lambda))}|r_{S,3-\ell_\delta^*}|^{-D_\delta^{\ell_\delta^*}} \right)}{-\log\delta} \\
&= \limsup_{\delta \to 0} D_\delta^{\ell_\delta^*}.
\end{align*}
Finally, the inequality above and \eqref{eq:2ub} imply that
\begin{equation}\label{eq:step2b}
	\dB(\Lambda)=\limsup_{\delta\to0}\max\{D_\delta^1,D_\delta^2\}.
\end{equation}

\subsection{Proof of the tertiary formulas in Theorem \ref{T1}} \label{subsec:tertform}
For $\ell \in \{1,2\}$, define $d_n^\ell$ with the equations
\begin{equation*}\begin{aligned}
\sum_{S \in \mathbb{G}_n} \left( \frac{|r_{S,\ell}|}{|r_{S,3-\ell}|} \right)^{\dB(\mathrm{p}_\ell\Lambda)} |r_{S,3-\ell}|^{d_n^\ell} = 1.
\end{aligned}\end{equation*}
Let us define $d_*^\ell := \limsup_{n \to \infty} d_n^\ell$ and $D_*^\ell:=\limsup_{\delta\to0}D_\delta^\ell$.

\begin{lemma} \label{lemma:bpn22}
For any $\eta \in [r_{\max},1)$,
\begin{align*}
D_*^\ell = D^\ell := &\inf \Bigg\{ \alpha > 0 \ \Bigg| \ \limsup_{k \to \infty} \Bigg( \sum_{S \in \mathbb{S}_{\eta^k}} \left( \frac{|r_{S,\ell}|}{|r_{S,3-\ell}|} \right)^{\dB(\mathrm{p}_\ell\Lambda)} |r_{S,3-\ell}|^{\alpha}\Bigg)^{1/k} < 1 \Bigg\}, \\
d_*^\ell = d^\ell := &\inf \Bigg\{ \alpha > 0 \ \Bigg| \ \limsup_{n \to \infty} \Bigg( \sum_{S \in \mathbb{G}_n} \left( \frac{|r_{S,\ell}|}{|r_{S,3-\ell}|} \right)^{\dB(\mathrm{p}_\ell\Lambda)} |r_{S,3-\ell}|^{\alpha}\Bigg)^{1/n} < 1 \Bigg\}.
\end{align*}
\end{lemma}
\begin{Proof}
For any $\alpha > D_*^\ell$ there exists $\varepsilon > 0$, $k_0 >0 $ such that for any $k > k_0$ satisfies $\alpha - \varepsilon > D_{\eta^k}^\ell$ (we will have the hierarchy: $\alpha>\alpha-\varepsilon>D_{\eta^k}^\ell$), and therefore by the definition of $D_{\eta^k}^\ell$, and using that $|r_{S,3-\ell}| \leq \eta^k$ for $S \in \mathbb{S}_{\eta^k}$, we will have that either
\begin{align*}
0=\limsup_{k \to \infty} \Bigg( \sum_{S \in \mathbb{S}_{\eta^k}} \left( \frac{|r_{S,\ell}|}{|r_{S,3-\ell}|} \right)^{\dB(\mathrm{p}_\ell\Lambda)} |r_{S,3-\ell}|^{\alpha}\Bigg)^{1/k} < 1
\end{align*} or
\begin{align*}
&0<\limsup_{k \to \infty} \Bigg( \sum_{S \in \mathbb{S}_{\eta^k}} \left( \frac{|r_{S,\ell}|}{|r_{S,3-\ell}|} \right)^{\dB(\mathrm{p}_\ell\Lambda)} |r_{S,3-\ell}|^{\alpha}\Bigg)^{1/k} \\
& \qquad < \eta^{-\varepsilon} \cdot \limsup_{k \to \infty} \Bigg( \sum_{S \in \mathbb{S}_{\eta^k}} \left( \frac{|r_{S,\ell}|}{|r_{S,3-\ell}|} \right)^{\dB(\mathrm{p}_\ell\Lambda)} |r_{S,3-\ell}|^{\alpha}\Bigg)^{1/k} \\
& \qquad = \limsup_{k \to \infty} \Bigg(\sum_{S \in \mathbb{S}_{\eta^k}} \left( \frac{|r_{S,\ell}|}{|r_{S,3-\ell}|} \right)^{\dB(\mathrm{p}_\ell\Lambda)} |r_{S,3-\ell}|^{\alpha} \eta^{-\varepsilon k}\Bigg)^{1/k} \\
& \qquad \leq \limsup_{k \to \infty} \Bigg( \sum_{S \in \mathbb{S}_{\eta^k}} \left( \frac{|r_{S,\ell}|}{|r_{S,3-\ell}|} \right)^{\dB(\mathrm{p}_\ell\Lambda)} |r_{S,3-\ell}|^{\alpha - \varepsilon}\Bigg)^{1/k} \\
& \qquad \leq \limsup_{k \to \infty} \Bigg( \sum_{S \in \mathbb{S}_{\eta^k}} \left( \frac{|r_{S,\ell}|}{|r_{S,3-\ell}|} \right)^{\dB(\mathrm{p}_\ell\Lambda)} |r_{S,3-\ell}|^{D_{\eta^k}^1}\Bigg)^{1/k} = 1.
\end{align*}
Thus $D^\ell \leq D_*^\ell$.

For the other direction, let $\alpha < D_*^\ell$, then for any $k_0 > 0$ there exists $k > k_0$ such that $\alpha< D_{\eta^k}^\ell$ and now
\begin{align*}
&\Bigg( \sum_{S \in \mathbb{S}_{\eta^k}} \left( \frac{|r_{S,\ell}|}{|r_{S,3-\ell}|} \right)^{\dB(\mathrm{p}_\ell\Lambda)} |r_{S,3-\ell}|^{\alpha}\Bigg)^{1/k} \geq \Bigg( \sum_{S \in \mathbb{S}_{\eta^k}} \left( \frac{|r_{S,\ell}|}{|r_{S,3-\ell}|} \right)^{\dB(\mathrm{p}_\ell\Lambda)} |r_{S,3-\ell}|^{D_{\eta^k}^1}\Bigg)^{1/k} = 1.
\end{align*}
This proves that $D^\ell \geq D_*^\ell$. For the proof of $d_*^\ell=d^\ell$, apply the same procedure.
\end{Proof}

\begin{lemma} \label{Lemma:bpn2} For $\ell \in \{1,2\}$, $d^\ell = D^\ell$.
\end{lemma}
\begin{Proof}
By Lemma \ref{lemma:bpn22}, it is enough to see that the two sums in their definitions differ only by a subexponential factor, and that is what we show. For the first direction, let $\eta \in [r_{\max},1)$ and $\iii \in \Delta_{\eta^k}$, then it is easy to see that
\begin{equation}\label{eq:bpn21}
k \cdot \frac{\log \eta}{\log r_{\min}} \leq |\iii| < k \cdot \frac{\log \eta}{\log r_{\max}} + 1.
\end{equation}
Hence, taking $\alpha > d^\ell$, we observe that
\begin{align*}
&\limsup_{k \to \infty} \Bigg( \sum_{S \in \mathbb{S}_{\eta^k}} \left( \frac{|r_{S,\ell}|}{|r_{S,3-\ell}|} \right)^{\dB(\mathrm{p}_\ell\Lambda)} |r_{S,3-\ell}|^{\alpha}\Bigg)^{1/k} \\
& \qquad \leq \limsup_{k \to \infty} \Bigg(\sum_{n = \big\lceil k \frac{\log \eta}{\log r_{\min}} \big\rceil }^{ \big\lceil k \frac{\log \eta}{\log r_{\max}} \big\rceil }\sum_{S \in \mathbb{G}_{n}} \left( \frac{|r_{S,\ell}|}{|r_{S,3-\ell}|} \right)^{\dB(\mathrm{p}_\ell\Lambda)} |r_{S,3-\ell}|^{\alpha}\Bigg)^{1/k} \\
& \qquad = \limsup_{k \to \infty} \Bigg(\max_{n \in\left\{ \big\lceil k \frac{\log \eta}{\log r_{\min}} \big\rceil,\ldots,\big\lceil k \frac{\log \eta}{\log r_{\max}} \big\rceil \right\} }\sum_{S \in \mathbb{G}_{n}} \left( \frac{|r_{S,\ell}|}{|r_{S,3-\ell}|} \right)^{\dB(\mathrm{p}_\ell\Lambda)} |r_{S,3-\ell}|^{\alpha}\Bigg)^{1/k} < 1.
\end{align*}
Thus $d^\ell \geq D^\ell$.

Absorbing \eqref{eq:bpn21} once again, we also have that for any $\iii\in\Sigma^{\{*\}}$, there exists at least one $k\in\N$ such that $S_{\iii}\in\mathbb{S}_{\eta^k}$. Any such $k$ must satisfy
\[ (|\iii|-1) \frac{\log r_{\max}}{\log \eta} < k \leq |\iii|\frac{\log r_{\min}}{\log \eta}.\]
Then proceed with the prequel argument to get $d^\ell \leq D^\ell$, ending the proof of Lemma \ref{Lemma:bpn2}. \end{Proof}

Lemma~\ref{Lemma:bpn2} with \eqref{eq:step2b} implies that
\begin{equation*}
	\dB(\Lambda)=\limsup_{n\to\infty}\max\{d_n^1,d_n^2\}.
\end{equation*}

\subsection{Proof of the quaternary formulas in Theorem \ref{T1}}

In the subsequent setup, we show that a single limit suffices for computation.

Let $s_n^\ell$ be the similarity dimension of $\mathrm{p}_\ell\mathbb{G}_n$. Since the \ref{AWSC} holds for the IFS $\mathrm{p}_\ell\mathbb{G}_n$, for any $\eps > 0$ there exists $N \in \N$ such that for all $n \geq N$
\[
s_n^\ell \in [\dB(\mathrm{p}_\ell\Lambda), \dB(\mathrm{p}_\ell\Lambda) + \eps].
\]
This can be seen by the definition of the \ref{AWSC} and by the proof of \cite[Theorem 4.2.16]{MR4661364}. Let $\mathfrak{D}_n^\ell$ be the unique solution of $1= \sum_{S \in \mathbb{G}_n} \left( \frac{|r_{S,\ell}|}{|r_{S,3-\ell}|} \right)^{s_n^\ell} |r_{S,3-\ell}|^{\mathfrak{D}_n^\ell}$. Now, 
\begin{align*}
1 &= \sum_{S \in \mathbb{G}_n} \left( \frac{|r_{S,\ell}|}{|r_{S,3-\ell}|} \right)^{\dB(\mathrm{p}_\ell\Lambda)} |r_{S,3-\ell}|^{d_n^{\ell}} \\
&\geq \sum_{S \in \mathbb{G}_n} \left( \frac{|r_{S,\ell}|}{|r_{S,3-\ell}|} \right)^{s_n^\ell} \min\Bigg\{ 1, \min_{\substack{S \in \mathbb{G}_n \\ |r_{S,\ell}| \geq |r_{S,3-\ell}|}}\left\{ \left(\frac{|r_{S,\ell}|}{|r_{S,3-\ell}|}\right)^{- \eps} \right\} \Bigg\} |r_{S,3-\ell}|^{d_n^\ell} \\
&\geq \sum_{S \in \mathbb{G}_n} \left( \frac{|r_{S,\ell}|}{|r_{S,3-\ell}|} \right)^{s_n^\ell} |r_{S,3-\ell}|^{d_n^\ell + \eps z},
\end{align*}
where
\begin{align*}
z := \max\left\{0, \Big(1-\frac{\log r_{\max}}{\log r_{\min}}\Big) \right\} \in [0,1].
\end{align*}
Since $\eps$ can be chosen arbitrarily small as $n \to \infty$, we can conclude that $\limsup_{n \to \infty} d_n^\ell \geq \limsup_{n \to \infty} \mathfrak{D}_n^\ell$. The proof for the other direction is analogous.

\vspace{15 pt}

\section{On the examples}
Here we derive the claimed dimension values for Example \ref{ex:1} and \ref{ex:2}. By Ngai and Wang \cite[Theorem~2.9]{NgaiWang}, and Nguyen \cite{Nguyen}, the principal projections satisfy the WSC for both examples.

\subsection{Example \ref{ex:1}}

By (4.62) from Bárány, Simon and Solomyak \cite{MR4661364}, for any $n \in \N^+$, 
\begin{equation*}
\Xi_n := \Big\{ \iii \in \{1,2,3\}^n \ \Big| \ \forall k \in \{1,\dots,n-1\}: \iii_k \iii_{k+1} \neq 13 \Big\}.
\end{equation*}
satisfies the definition of $\Gamma_2^{\{n\}}$. Let $\alpha = {\log_{4}3}$. For $\iii \in \Xi_n$, denote $(\#\{S\in\mathbb{G}_n\, \big|\, \mathrm{p}_{2} S= \mathrm{p}_2 S_{\iii}\})^{\alpha}$ by $R_{\iii}$. Thus
\begin{equation}\label{eq:andef}\begin{aligned}
\sum_{\iii\in\Gamma_2^{\{n\}}} (\#\{S\in\mathbb{G}_n\, \big|\, \mathrm{p}_{2} S = \mathrm{p}_2 S_{\iii}\})^{\alpha} &= \sum_{\iii \in\Gamma_2^{\{n\}}} R_{\iii} = \sum_{\iii \in \Xi_n} R_{\iii} \\&= \sum_{\substack{\ \ \iii \in \Xi_n\\ \iii_n = 1}} R_{\iii} + \sum_{\substack{\ \ \iii \in \Xi_n \\ \iii_n = 3}} R_{\iii} + \sum_{k=1}^{n-1}\sum_{\substack{\ \ \iii \in \Xi_n \\ \iii_n = 2 \\ \iii_{n-1} = 2 \\ \dots \\ \iii_{n-k+2} = 2 \\ \iii_{n-k+1} \neq 2}} R_{\iii} \\[-20pt]
& =: a_0^{(n)} + a_1^{(n)}  + \sum_{k=1}^{n-1} a_{k+1}^{(n)}
\end{aligned}\end{equation}
where the last line defined $a_k^{(n)}$,  $k \in \{ 0,1,\dots,n\}$ in order. Clearly, $a_i^{(n)} := 0$ for $i>n$. Denote $a^{(n)} := \big( a_0^{(n)}, a_1^{(n)}, a_2^{(n)}, \dots \big) \in \R^{\N}$. By Corollary \ref{T3},
\begin{equation*}
\dH(\Lambda) = \lim_{n \to \infty}\bigg\{\frac{1}{n}\log_{3}\| a^{(n)} \|_1\bigg\},
\end{equation*}
where $\|.\|_1$ is the usual 1-norm of real sequences.
Although the decomposition of $a^{(n)}$ may seem ad hoc, now we show what it represents: if $\iii \in \Xi_n$ ends with 1 or 3, then for any $j \in \{1,2,3\}$ we have $R_{\iii j} = R_{\iii}$. On the other hand, if $\iii$ ends with exactly $\ell$ 2s, then $R_{\iii 2} = R_{\iii 3} = R_{\iii}$, but $R_{\iii 1} = (\ell+1)^\alpha \cdot R_{\iii}$, since
\begin{equation*}
\mathrm{p}_2S_{\iii22\dots221} = \mathrm{p}_2S_{\iii22\dots213} = \mathrm{p}_2S_{\iii22\dots133} = \dots = \mathrm{p}_2S_{\iii13\dots333}
\end{equation*}
while $S_{\iii22\dots221}, S_{\iii22\dots213}, S_{\iii22\dots133}, \dots , S_{\iii13\dots333}$ are $\ell + 1$ different functions. Therefore,
\begin{equation*}\begin{aligned}
a_0^{(n)} = \sum_{\substack{\ \ \ \iii \in \Xi_n \\ \iii_n = 1}} R_{\iii} &= \sum_{\substack{\ \ \ \iii \in \Xi_n \\ \iii_n = 1 \\ \iii_{n-1} = 1}} R_{\iii} + \sum_{\substack{\ \ \ \iii \in \Xi_n \\ \iii_n = 1 \\ \iii_{n-1} = 3}} R_{\iii} + \sum_{k=1}^{n-1}\sum_{\substack{\ \ \ \iii \in \Xi_n \\ \iii_n = 1 \\ \iii_{n-1} = 2 \\ \dots \\ \iii_{n-k+1} = 2 \\ \iii_{n-k} \neq 2}} R_{\iii} \\
&= a_0^{(n-1)} + a_1^{(n-1)} + \sum_{k=1}^{n-1} (k+1)^\alpha a_{k+1}^{(n-1)} \\
&\qquad\qquad\qquad \Big( + 0 \text{ disguised as } \sum_{k=n}^{\infty} (k+1)^\alpha a_{k+1}^{(n-1)} \Big).
\end{aligned}\end{equation*}
Similarly,
\begin{equation*}\begin{aligned}
a_1^{(n)} = \sum_{\substack{\ \ \ \iii \in \Xi_n \\ \iii_n = 3}} R_{\iii} &= \sum_{\substack{\ \ \ \iii \in \Xi_n \\ \iii_n = 3 \\ \iii_{n-1} = 1}} R_{\iii} + \sum_{\substack{\ \ \ \iii \in \Xi_n \\ \iii_n = 3 \\ \iii_{n-1} = 3}} R_{\iii} + \sum_{k=1}^{n-1}\sum_{\substack{\ \ \ \iii \in \Xi_n \\ \iii_n = 3 \\ \iii_{n-1} = 2 \\ \dots \\ \iii_{n-k+1} = 2 \\ \iii_{n-k} \neq 2}} R_{\iii} \\
&= 0 + a_1^{(n-1)} + \sum_{k=1}^{n-1} a_{k+1}^{(n-1)} \Big(+ 0 \text{ disguised as } \sum_{k=n}^{\infty} a_{k+1}^{(n-1)} \Big).
\end{aligned}\end{equation*}
Moreover,
\begin{equation*}\begin{aligned}
a_2^{(n)} = \sum_{\substack{\ \ \ \iii \in \Xi_n \\ \iii_n = 2 \\ \iii_{n-1} \neq 2}} R_{\iii} = a_{0}^{(n-1)} + a_{1}^{(n-1)},
\end{aligned}\end{equation*}
while for $k \in \{3,\dots,n\}$:
\begin{equation*}\begin{aligned}
a_k^{(n)} = \sum_{\substack{\ \ \ \iii \in \Xi_n \\ \iii_n = 2 \\ \iii_{n-1} = 2 \\ \dots \\ \iii_{n-k+2} = 2 \\ \iii_{n-k+1} \neq 2}} R_{\iii} = \sum_{\substack{\ \ \ \iii \in \Xi_n \\ \iii_{n-1} = 2 \\ \iii_{n-2} = 2 \\ \dots \\ \iii_{n-k+2} = 2 \\ \iii_{n-k+1} \neq 2}} R_{\iii} = a_{k-1}^{(n-1)}.
\end{aligned}\end{equation*}
From these, we conclude that
\begin{equation*}
a^{(n)} = \mathbf{L}a^{(n-1)}= \dots = \mathbf{L}^n(1,1,1,0,0,\dots) = \mathbf{L}^{n+1}(0,1,0,\dots),
\end{equation*}
where we define the operator $\mathbf{L}: \R^{\N} \to \R^{\N}$ as
\begin{equation*}
\mathbf{L} = \begin{bmatrix}
1 & 1 & 2^\alpha & 3^\alpha & 4^\alpha & 5^\alpha & \cdots \\
0 & 1 & 1 & 1 & 1 & 1 & \cdots \\
1 & 1 & 0 & 0 & 0 & 0 & \cdots \\
0 & 0 & 1 & 0 & 0 & 0 & \cdots \\
0 & 0 & 0 & 1 & 0 & 0 & \cdots \\
0 & 0 & 0 & 0 & 1 & 0 & \cdots \\
\vdots & \vdots & \vdots & \vdots & \vdots & \vdots & \ddots
\end{bmatrix}.
\end{equation*}

\begin{lemma}\label{ex2lemma1}
There exists a unique $\lambda^*\in(1,\infty)$ such that there is $a \in \R^{\N}$ with positive entries, $a_0,a_1,a_2 \geq 1$ and with $\mathbf{L}a = \lambda^*a$. Furthermore
\begin{equation*}
\lambda^* = \frac{1}{\lambda^*-1}\sum_{k=2}^\infty k^\alpha (\lambda^*)^{2-k} + \frac{(\lambda^*)^2}{(\lambda^*-1)^3}.
\end{equation*}
\end{lemma}

\begin{Proof}
Let $a \in \R^{\N}$ be such that
\begin{equation*}
\sum_{k=2}^\infty k^\alpha a_k < \infty.
\end{equation*}
Suppose $\mathbf{L}a = \lambda a$, $\lambda \in (1,\infty)$, then
\begin{align} \label{eq34765}
\forall k \geq 3: \quad \lambda a_k = \mathbf{L}a_k = a_{k-1} \ \implies \ a_k = \lambda^{2-k}a_2
\end{align}
and
\begin{align*}
\lambda a_1 = \mathbf{L}a_1 = a_1 + \sum_{k=2}^\infty a_k = a_1 + \sum_{k=2}^\infty \lambda^{2-k}a_2 &= a_1 + a_2\frac{\lambda}{\lambda-1} \\
& \implies \ a_1 = a_2\frac{\lambda}{(\lambda-1)^2}.
\end{align*}
Then
\begin{equation*}\begin{aligned}
\lambda a_0 &= \mathbf{L}a_0 = a_0 + a_1 + \sum_{k=2}^\infty k^\alpha a_k = a_0 + a_2\frac{\lambda}{(\lambda-1)^2} + \sum_{k=2}^\infty k^\alpha \lambda^{2-k}a_2 \\
&\implies  a_0 = \frac{1}{\lambda-1}\Big(\frac{\lambda}{(\lambda-1)^2} + \sum_{k=2}^\infty k^\alpha \lambda^{2-k}\Big)a_2
\end{aligned}\end{equation*}
and
\begin{equation*}\begin{aligned}
\lambda a_2 = \mathbf{L}a_2 = a_0 + a_1, \quad \ a_1 = a_2\frac{\lambda}{(\lambda-1)^2}
\end{aligned}\end{equation*}
implies that
\begin{equation}\begin{aligned} \label{eq34526}
\lambda = \frac{1}{\lambda-1}\Big(\frac{\lambda}{(\lambda-1)^2} + \sum_{k=2}^\infty k^\alpha \lambda^{2-k}\Big) + \frac{\lambda}{(\lambda-1)^2} = \frac{1}{\lambda-1}\sum_{k=2}^\infty k^\alpha \lambda^{2-k} + \frac{\lambda^2}{(\lambda-1)^3}.
\end{aligned}\end{equation}
For $\lambda \in (1,\infty)$ the right-hand side of \eqref{eq34526} strictly decreases continuously. Thus \eqref{eq34526} is solved by a unique $\lambda^*$ on $(1,\infty)$. Hence, choosing $a_2$ sufficiently large, the claim on the existence of  $\mathbf{L}a = \lambda^*a$ also follows.
\end{Proof}

\begin{lemma}\label{ex2lemma2}
Let $a^{(n)}$ be the sequence defined in \eqref{eq:andef}. Then $\lim_{n\to\infty}\frac{1}{n}\log\|a^{(n)}\|_1 = \log \lambda^*$.
\end{lemma}

\begin{Proof}
	Let $M \in \N^+$, define $\mathbf{L}_M:\R^M \to \R^M$ as $(\mathbf{L} \circ \text{proj}_M)|_M$, where $\text{proj}_M$ is the projection of $\R^{\N}$ to the subspace spanned by the first $M$ coordinates. Then $\mathbf{L}_M$ can be represented as the non-negative, irreducible aperiodic, $M$ by $M$ matrix:
\begin{equation*}
\mathbf{L}_M = \begin{bmatrix}
1 & 1 & 2^\alpha & 3^\alpha & \cdots & (M-2)^\alpha & (M-1)^\alpha\\
0 & 1 & 1 & 1 & \cdots & 1 & 1 \\
1 & 1 & 0 & 0 & \cdots & 0 & 0 \\
0 & 0 & 1 & 0 & \cdots & 0 & 0 \\
\vdots & \vdots & \vdots & \vdots & \ddots & \vdots & \vdots\\
0 & 0 & 0 & 0 & \cdots & 0 & 0 \\
0 & 0 & 0 & 0 & \cdots & 1 & 0
\end{bmatrix}.
\end{equation*}
By the Perron-Frobenius Theorem, there exists a unique $\lambda_M > 0$ such that $\lim_{n \to \infty}\frac{1}{n}\log \|\mathbf{L}_M^nv\|_1 = \log \lambda_M$ for any $0 \neq v \in \R^M$ with non-negative entries, and there is a $v^* \in \R^M$ with positive entries such that $\mathbf{L}_Mv^* = \lambda_Mv^*$. Therefore, with computations similar to \eqref{eq34765}-\eqref{eq34526} we have:
\begin{equation*}\begin{aligned}
\forall k \in [3,M-1]: \quad \lambda_M v_k^* = v_{k-1}^* \ \implies \ v_k^* = \lambda_M^{2-k}v^*_2
\end{aligned}\end{equation*}
\begin{equation*}\begin{aligned}
\lambda_M v_1^* = v^*_1 + \sum_{k=2}^{M-1} v^*_k = v^*_1 + \sum_{k=2}^{M-1} \lambda_M^{2-k}v^*_2 \ \implies \ v^*_1 = v^*_2\frac{1}{\lambda_M-1}\sum_{k=2}^{M-1} \lambda_M^{2-k}.
\end{aligned}\end{equation*}
Then
\begin{align*}
\lambda_M v_0^* = v^*_0 + v^*_1 + \sum_{k=2}^{M-1} k^\alpha v^*_k \implies  v^*_0 = \Big(\frac{1}{(\lambda_M-1)^2}\sum_{k=2}^{M-1} \lambda_M^{2-k} + \frac{1}{\lambda_M-1}\sum_{k=2}^{M-1} k^\alpha \lambda_M^{2-k}\Big)v^*_2
\end{align*}
and
\begin{equation*}\begin{aligned}
\lambda_M v^*_2 = v^*_0 + v^*_1, \quad \ v^*_1 = v^*_2\frac{1}{\lambda_M-1}\sum_{k=2}^{M-1} \lambda_M^{2-k}
\end{aligned}\end{equation*}
implies that
\begin{equation}\begin{aligned} \label{eq324988z9}
\lambda_M = \frac{1}{v^*_2}(v^*_0 + v^*_1) = \frac{\lambda_M}{(\lambda_M-1)^2}\sum_{k=2}^{M-1} \lambda_M^{2-k} + \frac{1}{\lambda_M-1}\sum_{k=2}^{M-1} k^\alpha \lambda_M^{2-k}.
\end{aligned}\end{equation}
From \eqref{eq34526} and \eqref{eq324988z9} $\lim_{M \to \infty}\lambda_M = \lambda^*$. Then $\|\mathbf{L}^{n-1}a^{(1)} \|_1 \geq \|\mathbf{L}_M^{n-1}a^{(1)}|_M \|_1$ follows inductively on $n$, remembering that all entries of $a^{(1)}$ are non-negative. Whence
\begin{equation}\label{eq9732}
\liminf_{n \to \infty}\frac{1}{n} \log \|\mathbf{L}^{n-1}a^{(1)} \|_1 \geq \log \lambda_M \longrightarrow \log \lambda^*.
\end{equation}
Finally, let $a \in \R^{\N}$ be such that $\mathbf{L}a = \lambda^*a$, $a_0,a_1,a_2 \geq 1$ and such that $a_k > 0$ for $k \in \N$. Then
\begin{equation*}
\|\mathbf{L}^{n-1}a^{(1)} \|_1 \leq \|\mathbf{L}^{n-1}a \|_1 = (\lambda^*)^{n-1} \|a \|_1
\end{equation*}
hence
\begin{equation*}
\limsup_{n \to \infty}\frac{1}{n} \log \|\mathbf{L}^{n-1}a^{(1)} \|_1 \leq \log \lambda^*
\end{equation*}
which, along with \eqref{eq9732}, proves the statement.
\end{Proof}

Then $\dH(\Lambda)=\log_3\lambda^*$ follows by Corollary~\ref{T3} and Lemmas \ref{ex2lemma1} and \ref{ex2lemma2}.

\subsection{Example \ref{ex:2}}

Now, we turn to our second example. The argument has two main components:
\begin{itemize}
\item The defining IFS has exact overlaps on all levels of iterates after the first, but these are generated only by 2 equalities: $S_{14} = S_{21}$ and $S_{24} = S_{31}$.
\item $\mathrm{p}_2\Lambda$ has additional exact overlaps, generated by the equality $\mathrm{p}_2 S_{34} = \mathrm{p}_2 S_{41}$.
\end{itemize}
Although these are visually convincing, their proof will be given in Lemmas \ref{ex2l1} and \ref{ex2l2}.

\subsubsection{Box-counting dimension}

We first deduce the value of the box-counting dimension. To apply Corollary \ref{tbh}, it is enough to compute the quantities
\begin{align*}
 \frac{1}{n}\log\big(\#\mathrm{p}_2\mathbb{G}_n\big) \quad \text{and}\quad \frac{1}{n}\log\big(\#\mathbb{G}_n\big).
\end{align*}
For the growth rate of $\#\mathbb{G}_n$, we notice that the cylinders form two types of objects: a disjoint cylinder rectangle and 3 overlapping ones, forming a stair-shaped structure. We call these two shapes type 1 and type 2. Now, we can see that a type 1 produces, after one iteration, exactly one type 1 and one type 2, while a type 2 gives rise, in the next level, to a type 1 and 3 type 2s. Hence,
\begin{equation} \label{eq:matrix}
\#\mathbb{G}_n =
\begin{bmatrix}
1 & 0
\end{bmatrix}
\begin{bmatrix}
1 & 1 \\
1 & 3
\end{bmatrix}^n
\begin{bmatrix}
1\\
3
\end{bmatrix}.
\end{equation}
Thus, by the Perron-Frobenius Theorem, 
\[\lim_{n\to\infty}\frac{1}{n}\log \big(\#\mathbb{G}_n \big) = \log(\lambda),\]
where $\lambda$ is the largest eigenvalue, $2+\sqrt{2}$, of the matrix in \eqref{eq:matrix}. For $\mathrm{p}_2\mathbb{G}_n$, observe that for every $n\geq1$
\begin{equation*}
\#\mathrm{p}_2\mathbb{G}_n = 3\cdot \#\mathrm{p}_2(\mathbb{G}_{n-1}) + 1,\ \#\mathrm{p}_2(\mathbb{G}_0) = 1\quad 
\end{equation*}
and so, by induction $\#\mathrm{p}_2\mathbb{G}_n = \sum_{i = 0}^n3^i$. In particular, $\lim_{n\to\infty}\frac{1}{n}\log\big(\#\mathrm{p}_2\mathbb{G}_n\big)=\log3$ (which agrees with the observation that $\mathrm{p}_2\Lambda = [0,1]$, and hence $\dB(\mathrm{p}_2\Lambda) = 1$, while $r_2 = 1/3$).
We conclude that
\begin{equation*}\begin{aligned}
\dB(\Lambda) &= -\lim_{n \to \infty}\frac{1}{n}\Bigg\{\frac{\log\big(\#\mathrm{p}_2\mathbb{G}_n\big)}{\log|r_2|}\Big(1 - \frac{\log |r_2|}{\log |r_1|}\Big) + \frac{\log \big(\#\mathbb{G}_n\big)}{\log |r_1|}\Bigg\} \\
&= \frac{\log3}{\log3}\Big(1 - \frac{\log 3}{\log 4}\Big) + \frac{\log(2+\sqrt{2})}{\log 4} = 1.093295401221 \dots
\end{aligned}\end{equation*}

\subsubsection{Hausdorff dimension}
For the Hausdorff dimension, we wish to proceed similarly to the computation of Example \ref{ex:1}. As some parts are very similar, we emphasise only the differences. First, we give the symbolic description.

 Let $\Sigma = \{1,2,3,4\}$, $V := \{21,31\}$ and $W := \{ 21, 31, 41\}$.

\begin{lemma} \label{ex2l1}

For any $n \in \N^+$ let 
\begin{align*}
\Xi_n :=&\ \Big\{ \iii \in \Sigma^n \ \Big| \ \forall k \in \{1,\dots,n-1\}: (\iii_k \iii_{k+1}) \notin W \Big\} \\
\Theta_n :=&\ \Big\{ \iii \in \Sigma^n \ \Big| \ \forall k \in \{1,\dots,n-1\}: (\iii_k \iii_{k+1}) \notin V \Big\}.
\end{align*}
Now $\Xi_n$ satisfies the definition of $\Gamma_2^{\{n\}}$ and $\Theta_n$ satisfies the definition of $\Sigma^{\{n\}}$.
\end{lemma}

\begin{Proof}

We only prove that $\Xi_n$ satisfies the definition of $\Gamma_2^{\{n\}}$. The other assertion has very similar proof, so we omit it.

We will show that the map $f_n: \Xi_n \longrightarrow \mathrm{p}_2\mathbb{G}_n, \iii \mapsto \mathrm{p}_2S_{\iii}$ is injective and surjective. Recall that for $\iii = (\iii_1 \iii_2 \dots \iii_{|\iii|}) \in \Sigma^{|\iii|}, 0 \leq n < m \leq |\iii|$, we define
\[
\iii|_{(n,m]} := (\iii_{n+1} \iii_{n+2} \dots \iii_{m-1} \iii_m) \in \Sigma^{m-n}.
\]

We use induction to prove surjectivity and injectivity. The statement is straightforward for $n = 2$. Let $n \in \N^+$ and assume that $f_n$ is surjective and injective. First, we show that $f_{n+1}$ is surjective. Let $\iii \in \Sigma^{n+1}$ be arbitrary. Since $f_n$ is surjective, there is $\jjj \in \Xi_n$ such that $\mathrm{p}_2S_{\iii} = f_n(\jjj) \circ \mathrm{p}_2S_{\iii_{n+1}}$. If $(\jjj_n \iii_{n+1}) \notin W$, then the concatenation $\jjj\iii_{n+1}\in\Xi_{n+1}$ and so $f_{n+1}(\jjj\iii_{n+1})=f_n(\jjj)\circ \mathrm{p}_2S_{\iii_{n+1}}$. We now consider each remaining case (i.e., each element of $W$) separately. First, assume that $(\jjj_n \iii_{n+1}) = 31$, then 
\[ \mathrm{p}_2S_{\iii} = f_{n-1}(\jjj|_{(0,n-1]})\circ \mathrm{p}_2S_{31} = f_{n-1}(\jjj|_{(0,n-1]})\circ \mathrm{p}_2S_{24}\]
and $(\jjj|_{(0,n-1]}24) \in \Xi_{n+1}$ since $\jjj|_{(0,n-1]} \in \Xi_{n-1},\, (24) \notin W$, and since for any $(ij) \notin W$ with $j \neq 2$, $(\jjj_{n-1}2) \notin W$. 
The case when $(\jjj_n \iii_{n+1}) = (41)$ is similar.

Lastly, assume that $(\jjj_n \iii_{n+1}) = 21$, then $(\jjj \iii_{n+1}) = (\kkk 22\dots221)$, where $\kkk$ either the empty word or a finite word ending in $1$, $3$ or $4$. If $\kkk$ is the empty word or $\kkk_{|\kkk|}=1$, then
\[ \mathrm{p}_2S_{\iii} = \mathrm{p}_2S_{\jjj \iii_{n+1}} = \mathrm{p}_2S_{\kkk 22\dots221} = \mathrm{p}_2S_{\kkk 14\dots44}\]
and $\kkk 14\dots44 \in \Xi_{n+1}$. If $\kkk_{|\kkk|}=3$, then
\[ \mathrm{p}_2S_{\iii} = \mathrm{p}_2S_{\jjj \iii_{n+1}} = \mathrm{p}_2S_{\kkk_- 322\dots221} = \mathrm{p}_2S_{\kkk_- 314\dots44}=\mathrm{p}_2S_{\kkk_- 244\dots44}\]
and $\kkk_- 244\dots44 \in \Xi_{n+1}$. The case $\kkk_{|\kkk|}=4$ is similar.

Let $n \in \N^+$. For injectivity, we prove that if $\iii,\jjj \in \Sigma^n$ with $\iii \neq \jjj$ are such that $\mathrm{p}_2S_{\iii} = \mathrm{p}_2S_{\jjj}$, then $\iii \notin \Xi_n$ or $\jjj \notin \Xi_n$ (it might happen that neither of them is in $\Xi_n$). Then $\iii = (\iii \wedge \jjj) \iii'$ and $\jjj = (\iii \wedge \jjj) \jjj'$, where $|\iii'| = |\jjj'| \geq 2$, since $\mathrm{p}_2S_{\iii'} \neq \mathrm{p}_2S_{\jjj'}$ for $\iii' \neq \jjj' \in \Sigma$.
Since $\mathrm{p}_2 S_{\iii'|_{(0,2]}}((0,1)) \cap \mathrm{p}_2 S_{\jjj'|_{(0,2]}}((0,1)) \neq \emptyset$, we have that
\[ \{\iii'|_{(0,2]}, \jjj'|_{(0,2]} \} \in \big\{ \{13,21\}, \{14, 21\}, \{14,22\}, \{23,31\}, \{24, 31\}, \{24,32\}, \{33,41\}, \{34, 41\}, \{34,42\} \big\}. \]

Six of the above cases contain a forbidden word from $W$, thus we only need to consider the remaining three cases:
\[ 
\{\iii'|_{(0,2]}, \jjj'|_{(0,2]} \} \in \big\{\{14,22\},\{24,32\},\{34,42\} \big\} = \big\{ \{i4, (i+1)2\} \, \big| \, i \in \{1,2,3\}\big\}.
 \]
Fix $i\in\{1,2,3\}$ and assume
\[
\iii'|_{(0,2]} = (i4), \qquad \jjj'|_{(0,2]} = ((i+1)2).
\]
Since $\mathrm{p}_2S_{i4} \neq \mathrm{p}_2S_{(i+1)2}$, it follows that $|\iii'|=|\jjj'|\geq 3$.

For convenience, we write \( j^k = \underbrace{jj\dots jj}_k \) for $j\in\Sigma$ and $k\in\N$. Let $k\geq 2$ be the largest integer such that
\[
\iii'|_{(0,k]} = (i4^{k-1}), \qquad \jjj'|_{(0,k]} = ((i+1)2^{k-1}).
\] 
Since $\mathrm{p}_2 S_{i4^{k-1}} \neq \mathrm{p}_2 S_{(i+1)2^{k-1}}$, $k \neq |\iii'|$. Then, using that $\mathrm{p}_2S_{i4} = \mathrm{p}_2S_{(i+1)1}$, we have that
\[
\mathrm{p}_2S_{i4^{k-1}} = \mathrm{p}_2S_{(i+1)14^{k-2}} = \mathrm{p}_2S_{(i+1)214^{k-3}} = \dots = \mathrm{p}_2S_{(i+1)2^{k-2}1}.
\]
Thus, 
\[
\Big\{ (f,g) \in \Sigma^2 \,\Big|\, \mathrm{p}_2 S_{i4^{k-1}f}((0,1)) \cap \mathrm{p}_2 S_{(i+1)2^{k-1}g}((0,1)) \neq \emptyset\Big\} = \{(3,1),(4,1),(4,2)\}.
\]
By the maximality of $k$, $(f,g) \neq (4,2)$. The remaining two possibilities give the claim.

\end{Proof}

Let $\alpha := \log_43$, $\Theta_* := \cup_{n \in \N^+}\Theta_n$ and $\Xi_* := \cup_{n \in \N^+}\Xi_n$. By definition, $\mathrm{p}_2\colon \Theta_*\to\Xi_*$ is well defined mapping with
\[ \mathrm{p}_2 (\iii) = \mathrm{p}_2 (\jjj) \quad \text{ meaning that } \quad \mathrm{p}_2 S_{\iii} = \mathrm{p}_2 S_{\jjj}.\]

For $\iii \in \Xi_n$, denote $(\#\{S\in\mathbb{G}_n\, \big|\, \mathrm{p}_{2} S= \mathrm{p}_2 S_{\iii}\})^{\alpha}$ by $R_{\iii}$. In Example \ref{ex:1}, we decomposed $\sum_{\iii \in \Gamma_2^{\{n\}}}R_{\iii}$ with respect to its last symbols; now we decompose with respect to the first symbol. 


\begin{lemma} \label{ex2l2}
For $\iii \in \Xi_n$ and $j \in \Sigma$ such that $j\iii_1 \notin W$,
\begin{align*}
R_{(j \iii)} = \begin{cases}(\frac{k+2}{k+1})^\alpha \cdot R_{\iii}& \text{if } j = 3 \text{ and } \iii = 3^k4\kkk \text{ for } k \leq n-1,\, \kkk \in \Xi_{n-k-1}, \\ 
R_{\iii} & \text{else}. \end{cases}
\end{align*}
\end{lemma}

\begin{Proof}

By Lemma \ref{ex2l1}, for any $\iii \in \Xi_n$
\[ R_{\iii} = (\#\{S\in\mathbb{G}_n\, \big|\, \mathrm{p}_{2} S = \mathrm{p}_2 S_{\iii}\})^\alpha = (\#\{\jjj \in \Theta_n\, \big|\, \mathrm{p}_{2}(\jjj) = \iii\})^\alpha. \]

We make some observations first. 
Let us fix $j \in \{1,2\}$. Then we have the following two facts:
\begin{itemize}
\item For any $\ell \in \Sigma$ such that $j\ell \in \Theta_2$, $\{(ab)\in \Theta_2| \mathrm{p}_2(ab) = \mathrm{p}_2(j\ell) \} = \{j\ell\}$. Thus, for any $\iii \in \Xi_{n-1}$ such that $j\iii \in \Xi_n$
\begin{align} \label{eqex2ob1}
\big\{\jjj \in \Theta_{n}\, \big|\, \mathrm{p}_2(\jjj) = \mathrm{p}_2(j\iii)\big\} = \big\{j\jjj \in \Theta_{n}\, \big|\, \mathrm{p}_2(\jjj) = \mathrm{p}_2(\iii)\big\}.
\end{align}
In particular, for every $\iii\in\Xi_{n-1}$ such that $j\iii\in\Xi_n$
\[R_{j\iii} = R_{\iii}.\]
\item Let $j \neq k\in\Sigma$. If $\mathrm{p}_2(kc) = \mathrm{p}_2(ab)$ where $(kc), (ab) \in \Theta_2$, then $a \neq j$. Thus, if $k\iii, \jjj \in \Theta_n$ are such that $\mathrm{p}_2(k\iii) = \mathrm{p}_2(\jjj)$, then
\begin{align}\label{eqex2ob2}
\jjj_1 \neq j.
\end{align}
\end{itemize}

Let $\iii \in \Xi_{n-1}$ such that $4\iii \in \Xi_n$. Now, $\iii_1 \neq 1$ as $41 \in W$. Thus, by \eqref{eqex2ob2}
\[
\big\{\jjj \in \Theta_{n}\, \big|\, \mathrm{p}_2(\jjj) = \mathrm{p}_2(4\iii)\big\} = \big\{4\jjj \in \Theta_{n}\, \big|\, \mathrm{p}_2(\jjj) = \mathrm{p}_2(\iii)\big\}.
\]
Therefore,
\[R_{4\iii} = R_{\iii}.\]

Let $k\in \N, n\geq k$ and $\iii \in \Theta_{n-k}$ such that $\iii_1 = 2$. Then $3^k\iii \in \Theta_n$. Then, $3^k\iii = 3^k2\dots$ and since $\{\jjj \in \Theta_{k+1} | \mathrm{p}_2(\jjj) = \mathrm{p}_2(3^k2)\} = \{3^k2\}$ by \eqref{eqex2ob2}
\[R_{3^k\iii} = R_{\iii}.\]

Finally, let $\iii = 3^k4\kkk \in \Xi_n$, where $k \leq n-1$. Then one easily has $\{\jjj \in \Theta_{k+1}| \mathrm{p}_2(\jjj) = \mathrm{p}_2(3^k4)\} = \big\{3^\ell41^{k-\ell} \, \big|\, \ell \in \{0,1,\dots,k\}\big\}$ and since $3^k4\kkk \in \Xi_n$, either $\kkk$ is the empty word, or $\kkk_1 \neq 1$. Thus, by \eqref{eqex2ob2}
\[
\#\{\jjj \in \Theta_n\, \big|\, \mathrm{p}_{2}(\jjj) = \mathrm{p}_2 (3^k4\kkk)\} = \#\{\jjj \in \Theta_{k+1}\, \big|\, \mathrm{p}_{2}(\jjj) = \mathrm{p}_2 (3^k4)\} \cdot \#\{\jjj \in \Theta_{n-k-1}\, \big|\, \mathrm{p}_{2}(\jjj) = \mathrm{p}_2 (\kkk)\}.
\] 
Therefore, for every $k\geq 1$
\[R_{3^k4\iii} = (k+1)^\alpha \cdot R_{\iii}=\frac{(k+1)^\alpha}{k^\alpha} k^\alpha R_{\iii}=\frac{(k+1)^\alpha}{k^\alpha} R_{3^{k-1}4\iii}.\]

\end{Proof}

Let
\begin{align*}
\sum_{\iii \in \Gamma_2^{\{n\}}} R_{\iii} = \sum_{\iii \in \Xi_n} R_{\iii} &= \sum_{\substack{\ \iii \in \Xi_n \\ \iii_1 = 1}} R_{\iii} + \sum_{\substack{\ \iii \in \Xi_n \\ \iii_1 = 2}} R_{\iii} + \sum_{\substack{\ \iii \in \Xi_n \\ \iii_1 = 4}} R_{\iii} + \sum_{\text{rest}} R_{\iii} + \sum_{k=1}^{n-1} \sum_{\substack{\ \iii \in \Xi_n \\ \iii = 3^k4\kkk}} R_{\iii} \\
&=: a_1^{(n)} \quad \  + a_2^{(n)} \quad \ \ + a_4^{(n)} \quad \ \ \ + a_3^{(n)} \ \ \, + \sum_{k=1}^{n-1} a_{3,k}^{(n)},
\end{align*}
where we defined $a_\ell^{(n)}$ naturally. We denote $a^{(n)} := \big( a_1^{(n)}, a_2^{(n)}, a_4^{(n)}, a_3^{(n)},  a_{3,1}^{(n)}, a_{3,2}^{(n)},\dots \big) \in (\R^+)^{\N}$. By Lemma \ref{ex2l2},
\begin{equation*}
a^{(n)} = \mathbf{L}a^{(n-1)}= \dots = \mathbf{L}^n(1,1,1,1,0,\dots),
\end{equation*}
where 
\begin{equation*}
\mathbf{L} = \begin{bmatrix}
1 & 1 & 1 & 1 & 1 & 1 & 1 & 1 &\cdots \\
0 & 1 & 1 & 1 & 1 & 1 & 1 & 1 &\cdots \\
0 & 1 & 1 & 1 & 1 & 1 & 1 & 1 &\cdots \\
0 & 1 & 0 & 1 & 0 & 0 & 0 & 0 &\cdots \\
0 & 0 & 2^\alpha & 0 & 0 & 0 & 0 & 0 &\cdots \\
0 & 0 & 0 & 0 & (\tfrac{3}{2})^\alpha & 0 & 0 & 0 &\cdots \\
0 & 0 & 0 & 0 & 0 & (\tfrac{4}{3})^\alpha & 0 & 0 &\cdots \\
0 & 0 & 0 & 0 & 0 & 0 & (\tfrac{5}{4})^\alpha & 0 &\cdots \\
\vdots & \vdots & \vdots & \vdots & \vdots & \vdots & \vdots & \vdots &\ddots
\end{bmatrix}.
\end{equation*}

Then, the equation $\lambda a = \mathbf{L} a$ defines a system of equations for $\lambda$ and for the elements of $a$. Expressing $\lambda$, we obtain the equation
\begin{equation*}
\lambda^2 = 2\lambda + (\lambda-1)\sum_{k=0}^\infty (k+1)^\alpha \lambda^{-k},
\end{equation*}
if $\lambda > 1$. From this, one can compute $\lambda^*$ numerically.

\bibliographystyle{abbrv}
\bibliography{cikk}

\end{document}